\documentclass[12pt,a4paper]{article}
 \usepackage[top=40mm,bottom=40mm]{geometry}

\usepackage{amsmath,defns,euler,natbib,url,graphicx,microtype,reducecode}
\usepackage{amsfonts,mathtools}

\usepackage{hyperref} 
\hypersetup{pdftex,
            backref=true,
            hyperindex=true,
            colorlinks=true,
            citecolor=blue,
            bookmarks=true,
            breaklinks=true}

\newcommand{\uu}{{\bar u}}

\renewcommand{\vec}[1]{\text{\boldmath$#1$}}

\newcommand{\lp}{l} 
\ifx\de\undefined
\newcommand{\de}[2]{\mathchoice
  {\frac{d #2}{d #1}}
  {{d #2}/{d #1}}
  {{d #2}/{d #1}}
  {{d #2}/{d #1}}
  }
\fi
\newcounter{i}
\newcommand{\EE}{\ensuremath{\mathbb E}}
\DeclareMathOperator{\Span}{span}
\newcommand{\Eu}{\operatorname{Eu}}
\allowdisplaybreaks

\title{Multiscale modelling couples patches of two-layer thin fluid flow}

\author{Meng Cao\thanks{School of Mathematical Sciences,
University of Adelaide, South Australia 5005.  \protect\url{mailto:meng.cao@adelaide.edu.au} or \protect\url{mailto:mengcao1188216@gmail.com}}
\and 
A.~J. Roberts\thanks{School of Mathematical Sciences,
University of Adelaide, South Australia 5005.  \protect\url{mailto:anthony.roberts@adelaide.edu.au}}}

\begin{document}
    
\maketitle

\begin{abstract}
The multiscale gap-tooth scheme uses a given microscale simulator of complicated physical processes to enable macroscale simulations by computing only only small sparse patches.
This article develops the gap-tooth scheme to the case of nonlinear microscale simulations of thin fluid flow.
The microscale simulator is derived by artificially assuming the fluid film flow having two artificial layers but no distinguishing physical feature.
Centre manifold theory assures that there exists a slow manifold in the two-layer fluid film flow.
Eigenvalue analysis confirms the stability of the microscale simulator.
This article uses the gap-tooth scheme to simulate the two-layer fluid film flow.
Coupling conditions are developed by approximating the values at the edges of patches by neighbouring macroscale values.
Numerical eigenvalue analysis suggests that the gap-tooth scheme with the developed two-layer microscale simulator empowers feasible computation of large scale simulations of fluid film flows.
We also implement numerical simulations of the fluid film flow by the gap-tooth scheme.
Comparison between a gap-tooth simulation and a microscale simulation over the whole domain demonstrates that the gap-tooth scheme feasibly computes fluid film flow dynamics with computational savings.
\end{abstract}

\tableofcontents

\section{Introduction}
\label{sec:intro}

Mathematical equations describing geophysical fluid dynamics are typically written at the macroscale of kilometres. 
But the underlying turbulent flow and physics is best understood at the very much finer `microscale' sub-metre scale. 
We aim to empower scientists and engineers to use brief bursts of a given microscale simulator of wave-like dynamics on small patches of the space-time domain in order to make efficient and accurate macroscale simulations without ever knowing a macroscale closure.

Our modelling further develops the equation-free gap-tooth scheme \cite[e.g.]{Gear03, Samaey03b, Samaey:2005fk, Samaey2009} to empower novel simulation of wave-like systems over large time and space scales from a \emph{given} microscopic simulator.
Previously most multiscale modelling techniques have been developed for dissipative systems~\cite[e.g.]{E:2003fk, Kevrekidis:2003fk, Roberts:2005uq, Hou:2008fk}.
We suppose that the wave-like microscale simulator is computationally expensive so that only small time and spatial domain simulations are feasible: one example of future interest is direct numerical simulation of depth resolved turbulent fluid floods.
The aim is for the microscale simulator to provide the necessary data for the macroscopic computation, so that whenever the microscale simulator improves, then the overall macroscale simulation will correspondingly improve.

\cite{Cao2013} initiated the application of the gap-tooth scheme to  linear wave equations with weak dissipation.
\cite{Cao2014a} then developed and theoretically supported the gap-tooth method for more general dispersive and nonlinear wave-like systems, and as an indicative application and test, applied the methodology to a Smagorinski model of turbulent shallow water flow.
The generic key is to use polynomial interpolation of macroscale quantities, across the unsimulated gap between patches, to provide  coupling conditions on macroscale quantities on the edges of each patch of the microsimulator.
However, in most applications the microscale simulator will have many internal modes.
When there are such internal modes in the gap-tooth simulation, an outstanding issue is that at each time step we need to `lift' the macroscale coupling data to an appropriate microscale configuration.
As a first attempt to address this issue of lifting for wave-like dynamics, this paper uses the gap-tooth scheme to model viscous flow of a layer of fluid at moderate Reynolds number.

The flow of rainwater on the road, windscreen or other draining problems~\cite[e.g.]{Chang1987, Chang1994}, and paint and coating flows~\cite[e.g.]{Weinstein2004} are a few examples of fluid film flows.
Dynamics of such thin film flows have been studied extensively~\cite[e.g.]{Benjamin1957, Roberts1996, Roberts:1998fk, Roy:2002ys}.
The first aim of this article is to construct a two-layer model for such thin film flow as the microscale simulator.
The main reason to develop a two-layer model (Section~\ref{sec:micro}) is that it has microscale modes requiring lifting, but without the full complexity of fully resolved vertical structures.
However, the two-layer model is itself a novel model for fluid flows at moderate Reynolds number.

Consider a thin fluid flow of depth~$h(x,t)$ on an inclined plate with the slope~$\tan\theta$.
Section~\ref{sec:micro} artificially assumes two layers in the thin fluid flow, which have no distinguishing physical feature, as shown in Figure~\ref{twolayerDia}.
From the Navier--Stokes \pde{}s and the boundary conditions on the free surface and the flat substrate, Section~\ref{sec:theory} uses centre manifold theory~\cite[e.g.]{Roberts1988, Roberts2013a, Aulbach2000, Potzsche:2006uq} to derive a semi-slow two-layer model in the flow fields of depth~$h(x,t)$ and layer mean velocities~$\uu_1(x,t)$ in the lower layer and~$\uu_2(x,t)$ in the upper layer: nondimensionally the main parts of the model are
\begin{subequations}
\label{twolayer:eqS}
\begin{align}
\D th&={}-\frac12\left(\D{x}{h\uu_1}+\D{x}{h\uu_2}\right)\,,\label{twolayer:hfs}
\\
\D{t}{\uu_1}&\approx{}0.826\left(\tan\theta-\D xh\right)+\frac{1}{\re}\left(-19.3\frac{\uu_1}{h^2}+6.98\frac{\uu_2}{h^2}\right)
-1.48\uu_1\D{x}{\uu_1}
\,,\label{twolayer:u1fs}
\\
\D{t}{\uu_2}&\approx{}1.002\left(\tan\theta-\D xh\right)+\frac{1}{\re}\left(6.98\frac{\uu_1}{h^2}-5.36\frac{\uu_2}{h^2}\right)
-1.25\uu_1\D{x}{\uu_1}\,,\label{twolayer:u2fs}
\end{align}
\end{subequations}
for Reynolds number~$\re$, and the plate slope~$\tan\theta$. 
The right-hand sides of~\eqref{twolayer:eqS}, and the more refined version~\eqref{twolayer:fs}, include the effects of gravitational forcing, bed drag, nonlinear advection, and  dispersion.

The stability analysis in section~\ref{sec:eig} shows  instabilities at high wavenumber: consequently, section~\ref{sec:eig}  introduces an asymptotically consistent regularising operator to stabilise the two-layer model.
 
Section~\ref{sec:patch} applies the gap-tooth scheme to simulate the fluid film flow with the two-layer model being the microscale simulator. 
The coupling conditions are extended by the novel proposed lifting of the one-layer mean velocity~$\uu(x,t)$ to the two-layer  velocities~$\uu_1(x,t)$ and~$\uu_2(x,t)$.  
Comparisons between the gap-tooth simulation and the microscale simulation over the whole domain, Section~\ref{sec:num}, indicate that the gap-tooth scheme successfully simulates the dynamics of thin fluid flow at moderate Reynolds numbers.

\section{Construct the two-layer microscale simulator}
\label{sec:micro}

This section describes the derivation of the two-layer model.  First, section~\ref{sec:govern} lists the 2D continuity and Navier--Stokes equations of a layer thin fluid flow and the boundary conditions on the free surface and on the plate.
Second, section~\ref{sec:theory} embeds these equations in a family of equations with modified surface stress and mid-depth continuity so that an emergent two-layer slow manifold exists.
The computer algebra of Appendix~\ref{sec:app} then constructs the two-layer manifold model as summarised in section~\ref{sec:lowmodel}.

\subsection{Governing equations and boundary conditions for two layer thin film flow}
\label{sec:govern}

\begin{figure}
\begin{center}
\setlength{\unitlength}{0.7ex}
\begin{picture}(99,44)(9,4)
\put(5,0){\includegraphics[width=104\unitlength]{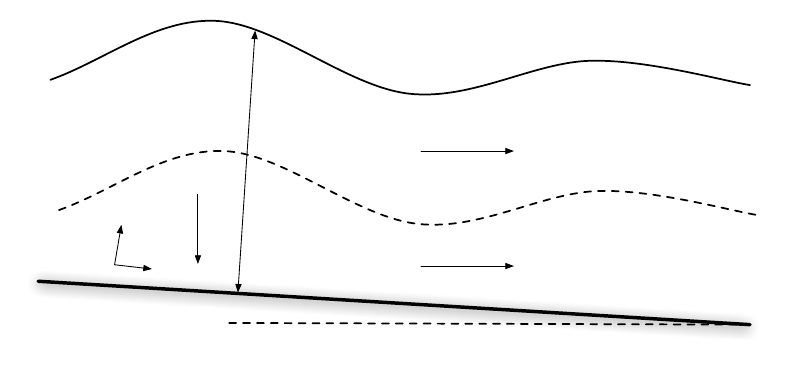}}
\put(25,12){$x$} \put(20,19){$z$}
\put(38,28){$h(x,z,t)$} \put(60,30){$\vec u_2(x,z,t)$} \put(60,15){$\vec u_1(x,z,t)$}
\put(32,15){$g$}
\put(90,7){$z=0$} \put(90,19){$z=h/2$}  \put(90,36){$z=h$}
\put(50,6){$\theta$}
\end{picture}
\end{center}
\caption{Diagram of the two-layer modelling of thin film flow.
Assume the fluid film flow has two artificial layers, which have no distinguishing physical feature.
Each layer has thickness~$h(x,t)/2$ and a mean layer velocity~$u_j(x,z,t)$ with $j=1$ for lower layer and $j=2$ for upper layer.
The plate has a slope~$\tan\theta$.
Then~$x$ and~$z$ are the lateral and normal coordinates, and~$g$ is gravity. }
\label{twolayerDia}
\end{figure}
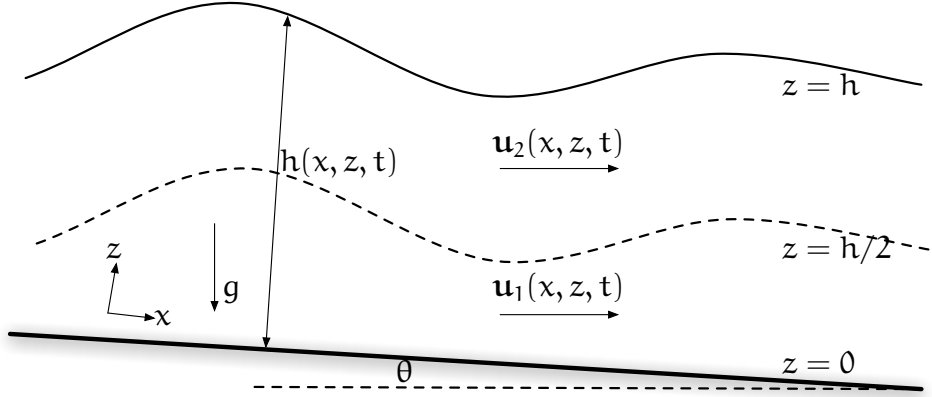

Consider a thin fluid flow of depth~$h(x,t)$ flowing down an inclined plate with slope~$\tan\theta$.
Denote the coordinate system by~$x$ and~$z$ along and normal to the plate. 
Assume the fluid film flow has two artificial layers, each of thickness~$h(x,t)/2$, which have no distinguishing physical feature, as shown schematically in Figure~\ref{twolayerDia}.
Each layer velocity field~$\vec q_j(x,z,t)=(u_j,w_j)$ and  pressure field by~$ p_j(x,z,t)$, where~$j=1$ is for lower layer and~$j=2$ is  for upper layer.

We nondimensionalise the system in terms of a typical film thickness~$H$, a typical fluid velocity $U:=\sqrt{gH}$, the constant density~\(\rho\), and the slope angle of the plate~$\theta$.
Then the Reynolds number $\re=HU/\nu_f$ with~$\nu_f$ being the viscosity of the fluid.
The nondimensional continuity and Navier--Stokes equations for the dynamics of the 2D fluid film flow are then
\begin{align}&
\divv\vec q_j=\D{x}{u_j}+\D{z}{w_j}=0\,,
\quad j=1,2\,,\label{twolayer:mass}
\\&
\D{t}{\vec q_j}+\vec q_j\cdot\grad\vec q_j=-\grad p_j+\frac{1}{\re}\grad^2\vec q_j+\vec g\,,
\quad j=1,2\,,\label{twolayer:mom}
\end{align}
where the Reynolds number~\(\re\) characterises the importance of the inertial terms compared to viscous dissipation, and where the vector $\vec g=(\tan\theta,-1)$ is the nondimensional forcing by gravity.

Well-known nondimensional boundary conditions hold on the plate and on the free surface.
\begin{itemize}
\item On the plate, prescribing no-slip requires
\begin{equation}
u_1=w_1=0
\quad\text{on }
z=0\,.
\label{twolayer:bcB}
\end{equation}
\item On the free surface, the pressure is assumed to be zero and so the stress normal to the surface is zero for zero surface tension \cite[e.g.]{Roberts1996, Roberts:1998fk}:
\begin{align}&
2\left[\D{z}{w_2}+\left(\D xh\right)^2\D{x}{u_2}-\D xh\left(\D{z}{u_2}+\D{x}{w_2}\right)\right]
\nonumber\\&\quad
-\re\left[1+\left(\D xh\right)^2\right] p_2=0
\quad\text{on }
z=h\,,\label{twolayer:bcNp}
\end{align}

\item The free surface having zero tangential stress results in
\begin{align}
\left[1-\left(\D xh\right)^2\right]\left(\D{z}{u_2}+\D{x}{w_2}\right)+2\left(\D xh\right)\left(\D{z}{w_2}-\D{x}{u_2}\right)=0
\quad\text{on } z=h\,.\label{twolayer:bcTpg}
\end{align}
\item Also on the free surface, the kinematic condition is
\begin{equation}
\D{t}{h}=w_2-u_2\D{x}{h}
\quad\text{on }
z=h\,.
\label{twolayer:bcK}
\end{equation}
\item Lastly, on the artificial interface of the two artificial layers, continuity of the physical fields requires
\begin{align}&
p_1=p_2\label{twolayer:bcI1}
\quad\text{on }
z=h/2\,,
\\& 
u_1=u_2\label{twolayer:bcI2}
\quad\text{on }
z=h/2\,,
\\&
w_1=w_2\label{twolayer:bcIw2}
\quad\text{on }
z=h/2\,,
\\&
\D{z}{u_1}=\D{z}{u_2}
\quad\text{on }
z=h/2\,.
\label{twolayer:bcI3g}
\end{align}
\end{itemize}
The \pde{}s~\eqref{twolayer:mass}--\eqref{twolayer:mom}, together with the boundary conditions \eqref{twolayer:bcB}--\eqref{twolayer:bcI3g} describe the dynamics of fluid film flow on an inclined plate.

\subsection{Embed to support with centre manifold theory}
\label{sec:theory}

This section embeds the system of physical equations in a family of artificial problems parametrised by~$\gamma$.
This embedding empowers theoretical support for the two-layer model.

First modify the tangential stress surface condition~\eqref{twolayer:bcTpg} to have an artificial forcing proportional to the net shear in the upper layer velocity:
\begin{align}&
\left[1-\left(\D xh\right)^2\right]\left(\D{z}{u_2}+\D{x}{w_2}\right)+2\left(\D xh\right)^2\left(\D{z}{w_2}-\D{x}{u_2}\right)
\nonumber\\&\quad
=(1-\gamma)\frac{2}{h}\left[u_2-u_2\left(x,{h}/{2},t\right)\right]
\quad\text{on }
z=h\,.\label{twolayer:bcTp}
\end{align}
When evaluated at parameter $\gamma=1$\,, the right-hand side of~\eqref{twolayer:bcTp} vanishes and so this artificial boundary condition~\eqref{twolayer:bcTp} reduces to the physical tangential stress conditions~\eqref{twolayer:bcTpg}.  
Also modify the derivative continuity~\eqref{twolayer:bcI3g} on the  interface of the two artificial layers to
\begin{align}
\left(1-\frac{\gamma}{2}\right)\D{z}{u_1}=\frac{\gamma}{2}\D{z}{u_2}+2(1-\gamma)\frac{u_1}{h}
\quad\text{on }
z=h/2\,.
\label{twolayer:bcI3}
\end{align}
When evaluated at $\gamma=1$, this artificial interface condition~\eqref{twolayer:bcI3} recovers the originally physical interface condition~\eqref{twolayer:bcI3g}. 
These are the two modifications of the embedding. 

The theoretical support is based upon a subspace of equilibria.
In the absence of lateral variations, \(\partial_x=0\)\,,  for parameter~$\gamma=0$\,, and on a horizontal bed, \(\theta=0\)\,, the fluid system~\eqref{twolayer:mass}--\eqref{twolayer:bcI3g} modified by equations~\eqref{twolayer:bcTp} and~\eqref{twolayer:bcI3} has two neutral modes of the layer shear flows $(u_1,u_2)\propto(z,h/2)$ and~$(u_1,u_2)\propto(0,z-h/2)$.
Conservation of fluid provides a third neutral mode in the dynamics.
That is, under the assumptions $\gamma=\partial_x=\theta=0$, there are three neutral modes, and hence a three dimensional `slow' subspace of equilibria, corresponding to uniform shear flows on a fluid of any thickness~$h$.
In the state space $(h,u_1,u_2,w_1,w_2,p_1,p_2)$ the slow subspace is
\begin{displaymath}
\EE_0=\Span\left\{
\begin{bsmallmatrix} 1\\0\\0\\0\\0\\0\\0 \end{bsmallmatrix},
\begin{bsmallmatrix}0\\z\\h/2\\0\\0\\0\\0  \end{bsmallmatrix}\,,
\begin{bsmallmatrix}0\\0\\z-h/2\\0\\0\\0\\0  \end{bsmallmatrix}
\right\}\,.
\end{displaymath}


Linearised about any equilibrium in the slow subspace~$\EE_0$ (together with the small parameters $\gamma=\partial_x=\theta=0$), equations~\eqref{twolayer:mass}, \eqref{twolayer:bcB} and~\eqref{twolayer:bcIw2} imply the normal velocity components $w_1=w_2=0$.
Then in the linearised system the normal component of the linearised \pde~\eqref{twolayer:mom} with \eqref{twolayer:bcNp} and~\eqref{twolayer:bcI1} establishes the pressures $p_1=p_2=0$.
Linearising the lateral component of the linearised \pde~\eqref{twolayer:mom} together with the modified boundary conditions\eqref{twolayer:bcI2} and~\eqref{twolayer:bcTp}--\eqref{twolayer:bcI3} gives a linear system for the lateral velocities~$u_j$:
\begin{subequations}
\label{eq:LSu}
\begin{align}&
\D{t}{u_j}-\frac{1}{\re}\DD{z}{u_j}=0\,,\label{twolayer:Lspec2}
\\&
u_1-u_2=0 \quad\text{on}\quad z=h/2\,,\label{twolayer:Lspec10}
\\&
\D{z}{u_2}-\frac{2}{h}\left[u_2-u_2|_{z=h/2}\right]=0 \quad\text{on}\quad z=h\,,\label{twolayer:Lspec11}
\\&
\D{z}{u_1}-\frac{2}{h}u_1=0\quad\text{on}\quad z=h/2\,.\label{twolayer:Lspec12}
\end{align}
\end{subequations}
The linearised \pde~\eqref{twolayer:Lspec2} implies the eigenvalue $\lambda=-k^2/\re$ by seeking the eigenvector of the lateral velocities $(u_1,u_2)\propto(0,\sin k(h-z))$. 
The boundary condition~\eqref{twolayer:Lspec11} imposes the constraint $kh/2=\tan kh/2$ with solutions $kh=0,8.986,15.451,\ldots$\,.
There exist a corresponding linearly independent generalised eigenvector with lower lateral velocity $u_1\neq0$.
A spectral gap exists between the three zero eigenvalues of the slow subspace~\(\EE_0\) and the non-zero eigenvalues headed by $\lambda\approx-80.763/(h^2\re)$.

Given the spectral gap, centre manifold theory for such `infinite dimensional' systems~\cite[e.g.]{Roberts1988, Roberts2013a, Aulbach2000, Potzsche:2006uq} supports the existence, emergence and construction of a slow manifold model based upon the slow subspace shear modes.
The slow manifold is constructed as a regular perturbation of the slow subspace \cite[]{Roberts1988, Potzsche:2006uq}.
Importantly, the theory supports the model in a finite domain of the parameters~\(\gamma\), \(\theta\) and~\(\partial_x\).
Evaluating the resulting slow manifold model at the physical case of artificial parameter \(\gamma=1\) then provides the model for the physical flow dynamics.
Table~\ref{twolater:CoefLmodel} [p.\pageref{twolater:CoefLmodel}] shows evidence that the modelling converges at \(\gamma=1\)\,.

\subsection{A low order model of the two layer flow}
\label{sec:lowmodel}

Computer algebra (Appendix~\ref{sec:app}) constructs the semi-slow manifold of the two layer thin fluid flow: we call it `semi-slow' because the model resolves two lateral shear modes, and because we  use `slow' to refer to a model that resolves only the gravest lateral shear mode.
The computer algebra program derives the semi-slow model in the flow fields of depth~$h(x,t)$, lower layer mean velocity~$\uu_1(x,t)$ and upper layer mean velocity~$\uu_2(x,t)$.

The order of errors in the construction is phrased in terms of small parameters.
In the theoretical support, the subspace of equilibria are found for \(\theta=\partial_x=\gamma=0\) and so these are necessarily small parameters of the semi-slow manifold.
Small lateral spatial derivative~$\partial_x$ corresponds to physical solutions varying slowly in~\(x\) and in such a context has recently been made a rigorous approximation \cite[]{Roberts2013a}.
Because the theoretical support is based upon the 3D subspace of equilibria  parametrised by~\(h\), \(\uu_1\) and~\(\uu_2\), the semi-slow model is formally global in~\(h\), \(\uu_1\) and~\(\uu_2\).
Nonetheless, we typically discard high order terms in~$\uu_j$ as being insignificant in practical parameter regimes \cite[e.g.]{Roberts:1998fk}.
Thus $\mathcal O(\uu_1^p+\uu_2^p+\partial_x^p+\theta^p)$ denotes the error terms for some exponent~$p$, which means each term explicitly expressed in the model has in total less than~$p$ factors  of these four parameters.
The bigger the exponent number~$p$, the higher the order of the modelling.
The artificial small parameter~$\gamma$ has no physical meaning but is crucial to rigorously establish the semi-slow manifold.
However, relatively high orders of the artificial parameter~$\gamma$ are required so that evaluating at~$\gamma=1$ is accurate.

Computer algebra in Appendix~\ref{sec:app} derives the physical flow fields of pressures~$p_1$ and~$p_2$, and layer velocities~$u_1$ and~$u_2$ in terms of the film thickness~$h$, layer mean velocities~$\uu_1$ and~$\uu_2$, and  scaled local normal coordinate $Z=z/h$:
\begin{subequations}\label{eqs:twolayerfld}%
\begin{align}
p_1&={}(1-Z)h
\nonumber\\&
{}+h^2\DD xh\left[-0.422+0.104Z-0.5Z^2)+\gamma(-0.0148-0.0107)\right]
\nonumber\\&
{}+ h\tan\theta\,\D{x}{h}\left[(0.0469-0.104Z)+\gamma(0.0375+0.0107Z)\right]
\nonumber\\&
{}+\frac{2}{\re}\left[\left(\D{x}{\uu_1}-\D{x}{\uu_2}-2\D{x}{\uu_1}Z\right)+\frac{1}{h}\D xh\left(\uu_1+\uu_2-2\uu_1Z\right)\right]
\nonumber\\&
{}+\gamma\frac{1}{\re}\left[(4.875\D{x}{\uu_1}-1.125\D{x}{\uu_2}-1.5\D{x}{\uu_1}Z+0.5\D{x}{\uu_2}Z)\right.
\nonumber\\&
{}\left.+\frac{1}{h}\D xh(-7.625\uu_1+2.375\uu_2+1.5\uu_1Z-0.5\uu_2Z)\right]
\nonumber\\&
{}+\mathcal O(\uu_1^3+\uu_2^3+\partial_x^3+\theta^3,\gamma^2)\,,\label{twolayer:p1Z}
\\
p_2&={}(1-Z)h
\nonumber\\&
{}+h^2\DD xh\left[(-0.417+0.115Z-0.5Z^2)+\gamma(-0.0178-0.0049Z)\right]
\nonumber\\&
{}+ h\tan\theta\,\D{x}{h}\left[(0.0521-0.115Z)+\gamma(0.0404+0.0049Z)\right]
\nonumber\\&
{}+\frac{4}{\re}\left[\left(-\D{x}{\uu_1}-\D{x}{\uu_2}+2\D{x}{\uu_1}Z\right)+\frac{1}{h}\D xh\left(2\uu_1-\uu_2-2\uu_1Z+\uu_2Z\right)\right]
\nonumber\\&
{}+\gamma\frac{1}{\re}\left[(11\D{x}{\uu_1}-3\D{x}{\uu_2}-13.75\D{x}{\uu_1}Z+4.25\D{x}{\uu_2}Z)\right.
\nonumber\\&
{}\left.+\frac{1}{h}\D xh(-13.75\uu_1+4.25\uu_2+13.75\uu_1Z-4.25\uu_2Z)\right]
\nonumber\\&
{}+\mathcal O(\uu_1^3+\uu_2^3+\partial_x^3+\theta^3,\gamma^2)\,,\label{twolayer:p2Z}
\\
u_1(Z)&={}\left[(4.0+1.5\gamma)\uu_1Z-0.5\gamma\uu_2Z+(4.0\uu_2-12.0\uu_1)\gamma Z^3\right]
\nonumber\\&
{}+\re h^2\D xh\left[(-0.104+0.0107\gamma)Z+0.5Z^2\right.
\nonumber\\&
{}\left.-(0.5+0.123\gamma)Z^3+0.225\gamma Z^5\right]
\nonumber\\&
{}+\re h^2 \tan\theta\,\left[(0.104-0.0107\gamma)Z-0.5Z^2\right.
\nonumber\\&
{}\left.+(0.5+0.123\gamma)Z^3-0.225\gamma Z^5\right]
\nonumber\\&
{}+\mathcal O(\uu_1^3+\uu_2^3+\partial_x^3+\theta^3,\gamma^2)\,,\label{twolayer:u1Z}
\\
u_2(Z)&={}\left[(6.0-11.9\gamma)\uu_1+(-2.0+4.13\gamma)\uu_2+(4.0-17.8\gamma)\uu_2Z
\right.
\nonumber\\&
{}\left.+(-8.0+48.3\gamma)\uu_1Z+(27.0\uu_2-69.0\uu_1)\gamma Z^2\right.
\nonumber\\&
{}\left.+(34.0\uu_1-14.0\uu_2)\gamma Z^3\right]
\nonumber\\&
{}+\re h^2\D xh\left[(-0.089+0.10\gamma)+(0.45-0.68\gamma)Z\right.
\nonumber\\&
{}\left.+(-0.63+1.8\gamma)Z^2+(0.25-2.5\gamma)Z^3+1.8\gamma Z^4-0.49\gamma Z^5\right]
\nonumber\\&
{}+\re h^2 \tan\theta\,\left[(0.089-0.10\gamma)+(-0.45+0.68\gamma)Z\right.
\nonumber\\&
\left.+(0.63-1.8\gamma)Z^2+(-0.25+2.5\gamma)Z^3-1.8\gamma Z^4+0.49\gamma Z^5\right]
\nonumber\\&
{}+\mathcal O(\uu_1^3+\uu_2^3+\partial_x^3+\theta^3,\gamma^2)\,.\label{twolayer:u2Z}
\end{align}
\end{subequations}
Equations~\eqref{eqs:twolayerfld} describe the low order shape of the manifold in the state space.
Physically, upon setting the parameter~$\gamma=1$, these four equations~\eqref{eqs:twolayerfld} approximately describe the vertical structures of the pressures and layer velocities associated with the terms of depth~$h(x,t)$, layer mean velocities~$\uu_j(x,t)$ and their lateral derivatives.
One important feature of this approach is that we do not impose these vertical structures on the flow: instead we systematically solve the governing physical fluid equations to discover the  structure functions appropriate for any suitable parameter regime.

The computer algebra in Appendix~\ref{sec:app} also determines the evolutions of the depth~$h(x,y)$ and layer mean velocities~$\uu_1(x,t)$ and~$\uu_2(x,t)$ on this semi-slow manifold~\eqref{eqs:twolayerfld} but now to high order in the parameter~$\gamma$:
\begin{subequations}\label{eqs:twolayercmm}%
\begin{align}
\D th&={}-0.5\D{x}{}\left(h\uu_1+h\uu_2\right)\,,\label{twolayer:hL}
\\
\D{t}{\uu_1}&={}\left(0.75+0.0438\gamma+0.0365\gamma^2-0.00439\gamma^3+0.0000522\gamma^4\right.
\nonumber\\&
{}\left.-0.000305\gamma^5-0.0000393\gamma^6\right)\left(\tan\theta-\D xh\right)
\nonumber\\&
{}+\frac{1}{\re}\left[-\left(18\gamma+1.35\gamma^2+0.0723\gamma^3-0.0869\gamma^4-0.0112\gamma^5-0.00637\gamma^6\right)\frac{\uu_1}{h^2}\right.
\nonumber\\&
{}\left.+\left(6.0\gamma+1.05\gamma^2-0.038\gamma^3-0.022\gamma^4-0.008\gamma^5-0.00218\gamma^6\right)\frac{\uu_2}{h^2}\right]
\nonumber\\&
{}+\mathcal O(\uu_1^3+\uu_2^3+\partial_x^3+\theta^3,\gamma^7)\,,\label{twolayer:u1L}
\\
\D{t}{\uu_2}&={}\left(1.125-0.195\gamma+0.0740\gamma^2-0.00126\gamma^3-0.0003062\gamma^4\right.
\nonumber\\&
{}\left.-0.000185\gamma^5-0.0000231\gamma^6\right)\left(\tan\theta-\D xh\right)
\nonumber\\&
{}+\frac{1}{\re}\left[\left(15\gamma-8.29\gamma^2+0.22\gamma^3+0.0278\gamma^4+0.0135\gamma^5+0.00351\gamma^6\right)\frac{\uu_1}{h^2}\right.
\nonumber\\&
{}\left.+\left(-9.0\gamma+3.17\gamma^2-0.0456\gamma^3-0.0173\gamma^4-0.00559\gamma^5-0.00146\gamma^6\right)\frac{\uu_2}{h^2}\right]
\nonumber\\&
{}+\mathcal O(\uu_1^3+\uu_2^3+\partial_x^3+\theta^3,\gamma^7)\,,\label{twolayer:u2L}
\end{align}
\end{subequations}
Equation~\eqref{twolayer:hL} is a direct consequence of conservation of fluid.
The momentum equations~\eqref{twolayer:u1L}--\eqref{twolayer:u2L} include the effects of viscous drag and dissipation~$\uu_j/h^2$, and the gravity forcing~$\tan\theta-\D xh$.
The higher order physical effects of advection and dispersion effects are included in the next equations~\eqref{twolayer:u1f}--\eqref{twolayer:u2f} where we report on the model to errors $\mathcal O(\uu_1^4+\uu_2^4+\partial_x^4+\theta^4,\gamma^7)$.

\begin{table}
\caption{Partial sums from evaluating coefficients at~$\gamma=1$ of selected terms in equation~\eqref{twolayer:u1L}--\eqref{twolayer:u2L} indicates that the power series in~$\gamma$ converges quickly.}
\label{twolater:CoefLmodel}
\def\o{\phantom{0}}
\begin{displaymath}
\begin{array}{ccccccc}
\hline
 & \parbox{0.1\linewidth}{$\D xh$ in~\eqref{twolayer:u1L}} 
 &  \parbox{0.1\linewidth}{ $\D xh$ in~\eqref{twolayer:u2L}}  
 &  \parbox{0.13\linewidth}{ $\uu_1/h^2/\re$ in~\eqref{twolayer:u1L}}  
 & \parbox{0.13\linewidth}{ $\uu_1/h^2/\re$ in~\eqref{twolayer:u2L}}
 & \parbox{0.13\linewidth}{ $\uu_2/h^2/\re$ in~\eqref{twolayer:u1L}} 
 & \parbox{0.13\linewidth}{$\uu_2/h^2/\re$ in~\eqref{twolayer:u2L}}
 \\
\hline
 \gamma^0 & -0.75\o\o & -1.125\, & 0 & 0 & 0 & 0
 \\
\gamma^1 & -0.7938 & -0.930& -18.\o\o & 15.\o\o & 6.\o\o\o & -9.\o\o\o
\\
\gamma^2 & -0.8302 & -1.004 & -19.35 & 6.712 & 7.050 & -5.288
\\
\gamma^3 & -0.8258 & -1.002 & -19.42 & 6.933 & 7.012 & -5.333
\\
\gamma^4 & -0.8259 & -1.002 & -19.34 & 6.960 & 6.990 & -5.350
\\
\gamma^5 & -0.8256 & -1.002 & -19.32 & 6.974 & 6.982 & -5.356
\\
\gamma^6 & -0.8255 & -1.002 & -19.32 & 6.977 & 6.980 & -5.357
\\
\hline
\end{array}
\end{displaymath}
\end{table}

Equations~\eqref{twolayer:p1Z}--\eqref{twolayer:u2L} express the semi-slow manifold model in terms of the introduced artificial parameter~$\gamma$.
Every coefficient in these equations is a power series in~$\gamma$. 
The partial sums in Table~\ref{twolater:CoefLmodel} indicate that these coefficient series in~$\gamma$ converges quickly for~$\gamma=1$.
Thus, when the model is constructed to errors~$\mathcal O(\gamma^7)$, it is apparent that all the shown digits are accurate.
\cite{Roberts1996,Roberts:1998fk} and~\cite{Roberts2002} reported similar convergence in other related physical problems.
Hereafter we calculate every coefficient in the model up to errors~$\mathcal O(\gamma^7)$, and then evaluate at the artificial parameter~$\gamma=1$.

Truncating to errors~$\mathcal O(\uu_1^4+\uu_2^4+\partial_x^4+\theta^4,\gamma^7)$, omitting the intricate details of the derivation and upon setting the artificial parameter~$\gamma=1$, the  evolution of the depth~$h(x,t)$, the lower layer mean velocity~$\uu_1(x,t)$ and the upper layer mean velocity~$\uu_2(x,t)$ on the semi-slow manifold are described by the flow conservation equation and by effective lateral momentum equations
\begin{subequations}
\label{twolayer:fs}
\begin{align}
\D th&={}-0.5\left(\D{x}{h\uu_1}+\D{x}{h\uu_2}\right)\,,\label{twolayer:hf}
\\
\D{t}{\uu_1}&\approx{}0.826\left(\tan\theta-\D xh\right)+\frac{1}{\re}\left(-19.3\frac{\uu_1}{h^2}+6.98\frac{\uu_2}{h^2}\right)
\nonumber\\&
{}-1.48\uu_1\D{x}{\uu_1}-0.225\uu_2\D{x}{\uu_2}+0.142\uu_2\D{x}{\uu_1}+0.0728\uu_1\D{x}{\uu_2}
\nonumber\\&
{}+\frac{(\uu_1-\uu_2)}{h}\left(-0.25\uu_1+0.34\uu_2\right)\D xh
\nonumber\\&
{}+\frac{1}{\re}\left(-3.84\DD{x}{\uu_1}+2.52\DD{x}{\uu_2}\right)\,,\label{twolayer:u1f}
\\
\D{t}{\uu_2}&\approx{}1.002\left(\tan\theta-\D xh\right)+\frac{1}{\re}\left(6.98\frac{\uu_1}{h^2}-5.36\frac{\uu_2}{h^2}\right)
\nonumber\\&
{}-1.25\uu_1\D{x}{\uu_1}-1.57\uu_2\D{x}{\uu_2}+0.768\uu_2\D{x}{\uu_1}+0.930\uu_1\D{x}{\uu_2}
\nonumber\\&
{}+\frac{(\uu_1-\uu_2)}{h}\D xh\left(-0.78\uu_1+0.38\uu_2\right)
\nonumber\\&
{}+\frac{1}{\re}\left(-1.98\DD{x}{\uu_1}+5.23\DD{x}{\uu_2}\right)\,.\label{twolayer:u2f}
\end{align}
\end{subequations}
Equation~\eqref{twolayer:hf} is a direct consequence of conservation of fluid.
The momentum equations~\eqref{twolayer:u1f}--\eqref{twolayer:u2f} include the effects of gravity~$\tan\theta-\D xh$, viscous drag~$\uu_j/h^2$, advection~$\uu_j\D{x}{\uu_i}$, dispersion~$\DD{x}{\uu_j}$, and other viscous terms, such as~$\D xh\D{x}{\uu_j}$.
Compared with the one-layer models~\cite[e.g.]{Prokopiou:1991fk, Roberts1996, Ruyer-Quil1998}, the two-layer model~\eqref{twolayer:u1f}--\eqref{twolayer:u2f} ensures more subtle effects and resolves more internal modes. 
These internal modes are necessary for a more generic microscale simulation in the gap-tooth scheme over that initiated for wave-like systems by \cite{Cao2014a}.


\subsection{Eigenanalysis of the microscale model}
\label{sec:eig}

This section linearly analyses the two-layer model~\eqref{twolayer:fs}.
Linear analysis indicates that an unphysical instability appears for high wavenumber.
\cite{Cao2014} explored three methods to avoid such instability: only resolving low wavenumbers; adding high order dissipation terms; and introducing a regularising operator.
We implemented, recommend and describe the last.

\paragraph{Linear analysis of the two-layer model}

Consider the modelled fluid film flow with two artificial layers on a flat plate with slope~$\tan\theta$.
The fluid shear flow has an equilibrium with thickness $h=1$\,, without loss of generality.
The model~\eqref{twolayer:fs}, assuming~$\partial_t=\partial_x=0$, simplifies to
\begin{align*}&
0.826\tan\theta+\frac{1}{\re}\left(-19.3\uu_1+6.98\uu_2\right)=0\,,
\\&
1.002\tan\theta+\frac{1}{\re}\left(6.98\uu_1-5.36\uu_2\right)=0\,,
\end{align*}
which predicts an equilibrium of layer mean velocities
\begin{align}
\uu_1=0.209\re\tan\theta
\quad\text{and}\quad
\uu_2=0.459\re\tan\theta\,.
\label{twolayer:eigEq}
\end{align}
Impose small perturbations to this equilibrium and seek solutions in the form
\begin{align}&
h=1+\hat h\exp(\lambda t+ikx)\,,
\nonumber\\&
\uu_1=0.209\re\tan\theta+\hat u_1\exp(\lambda t+ikx)\,,\label{twolayer:eigSol}
\\&
\uu_2=0.459\re\tan\theta+\hat u_2\exp(\lambda t+ikx)\,,\nonumber
\end{align}
for growth rate~$\lambda$~(possibly complex) and nondimensional wavenumber~$k$.

Substitute the form~\eqref{twolayer:eigSol} into the model~\eqref{twolayer:fs}, equate coefficients and derive the linear problem
\begin{align*}
\lambda\begin{bmatrix}\hat h \\ \hat u_1\\ \hat u_2 \end{bmatrix}=M\begin{bmatrix}\hat h \\ \hat u_1\\ \hat u_2 \end{bmatrix}\,,
\end{align*}
where the coefficient matrix~$M$ is
\begin{align}
\begin{bmatrix} 
\parbox{0.28\linewidth}{\raggedright$-0.334\tan\theta\re ik$}&
\parbox{0.28\linewidth}{\raggedright$ -0.5ik $}&
\parbox{0.28\linewidth}{\raggedright$ -0.5ik $}
\\
\\
\parbox{0.28\linewidth}{\raggedright$-(0.826+0.0259\tan\theta\re)ik$}&
\parbox{0.28\linewidth}{\raggedright $ \frac{1}{\re}(-19.3+3.84k^2)-0.244\tan\theta\re ik$}&
\parbox{0.28\linewidth}{\raggedright $ \frac{1}{\re}(6.98-2.52k^2)-0.088\tan\theta\re ik$}
\\
\\
\parbox{0.28\linewidth}{\raggedright$-(1.002+0.0029\tan\theta\re)ik$}&
\parbox{0.28\linewidth}{\raggedright $ \frac{1}{\re}(6.98+1.98k^2)+0.935\tan\theta\re ik$}&
\parbox{0.28\linewidth}{\raggedright $ \frac{1}{\re}(-5.36-5.23k^2)-0.526\tan\theta\re ik$}
 \end{bmatrix}.
\label{twolayer:eigmatric}
\end{align}
The coefficient matrix~$M$ has characteristic equation
\begin{align}&
\lambda^3+\lambda^2\left[1.114\tan\theta\re ik+\frac{1}{\re}\left(1.39k^2+24.66\right)\right]
\nonumber\\&\quad
-\lambda\left[0.476\tan^2\theta\re^2k^2+\tan\theta k(-2.265ik^2-14.03i-0.0144\re k)\right.
\nonumber\\&\quad
\left.+15.09\frac{1}{\re^2}k^4-0.914k^2-84.13\frac{1}{\re^2}k^2-54.73\frac{1}{\re^2}\right]
\nonumber\\&\quad
-\left[0.0725\tan^3\theta\re^3ik^3+\tan^2\theta\re k^2(-0.0191\re ik+0.615k^2+1.908)\right.
\nonumber\\&\quad
\left.+\tan\theta\frac{1}{\re}k(5.19ik^4-0.682\re^2ik^2-28.94ik^2-18.83i-0.841\re k^3\right.
\nonumber\\&\quad
\left.-0.198\re k)+\frac{1}{\re}k^2(0.209k^2-18.26)\right]=0\,.
\label{twolayer:eigfun}
\end{align}

\begin{figure}
\centering
\begin{tabular}{c@{\ }c}
\rotatebox{90}{\hspace{15ex}growth rate~$\Re\lambda$} &
\includegraphics{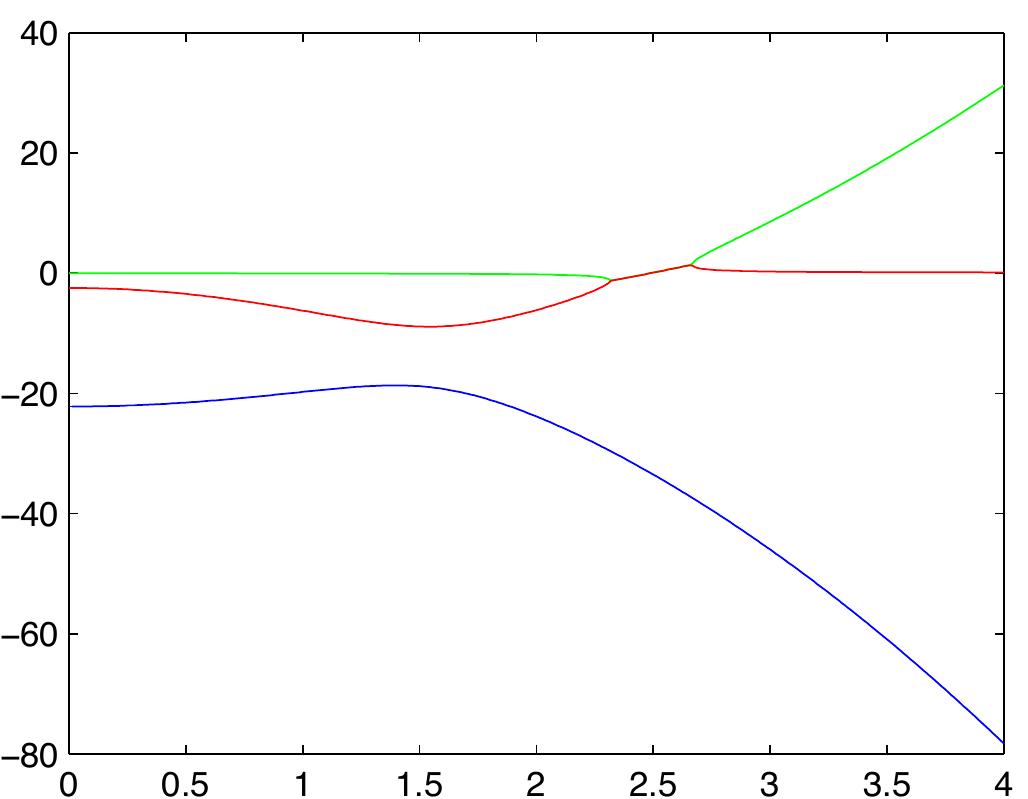}\\
& wavenumber~$k$
\end{tabular}
\caption{The growth rate~$\Re\lambda$ varies with the nondimensional wavenumber~$k$ from the characteristic equation~\eqref{twolayer:eigfun} of the model~\eqref{twolayer:fs}.
The slope of the plate is $\tan\theta=0$.
The Reynolds number is $\re=1$.
An unphysical instability appears for high nondimensional wavenumber $k>2.5$.
The little interval near $k=2.5$ is due to the numerical error.}
\label{twolayerEig1}
\end{figure}

Figure~\ref{twolayerEig1} plots the growth rates~$\Re\lambda$ versus the nondimensional wavenumber~$k$ from equation~\eqref{twolayer:eigfun}. 
All the values represented by the blue curve are negative which nicely reflects viscous decay of lateral shear modes. 
When the nondimensional wavenumber~$k<2.5$, the values represented by the red curve are negative~(viscous decay), and by the green curve are zeros~(conservation of fluid). 
But when the nondimensional wavenumber~$k>2.5$, the green curve increases to positive which implies that instability arises in the system.
This instability is nothing to do with physical instabilities, for example, found by~\cite{Chen:1993fk}, who found the instability of the two liquid films down an inclined plate due to the different viscosity at the interface of the two layer flow and on the free surface.
Here the instability arises at high wavenumber, a wavelength comparable to the thickness of the fluid, whereas our modelling is accurate for low wavenumber, long length scales.


\paragraph{Consistently avoid the instability}

We introduce a regularising operator to stabilise the unphysical instability.
For the momentum equations~\eqref{twolayer:u1f}--\eqref{twolayer:u2f}, consider applying the regularising operator $\mathcal L=1-C\partial_x(h^2\partial_x)$ to both sides of both \pde{}s.
The coefficient~$C$ is positive. 
The reason for using~$h^2$ in the regularising operator~$\mathcal L$ is to be dimensionally consistent which means we can cancel the~$h^{-2}$ in the drag terms in equations \eqref{twolayer:u1f}--\eqref{twolayer:u2f} and also ensure the regularising operator~$\mathcal L$ is self-adjoint.
This regularising operator  generates dissipation effects to counteract the positive growth rates.
The computer algebra of Appendix~\ref{sec:app} applies the regularising operator to the momentum \pde{}s~\eqref{twolayer:u1f}--\eqref{twolayer:u2f}, and gives the \pde{}s, in term of the regularising parameter~\(C\),
\begin{subequations}\label{eqs:twolayerI}%
\begin{align}
\D th&={}-0.5\left(\D{x}{h\uu_1}+\D{x}{h\uu_2}\right)\,,\label{twolayer:hI}
\\
\mathcal L\D{t}{\uu_1}&\approx{}0.826\left(\tan\theta-\D xh\right)
\nonumber\\&
{}+\frac{1}{\re}\left(-19.3\frac{\uu_1}{h^2}+6.98\frac{\uu_2}{h^2}\right)
\nonumber\\&
{}-1.48\uu_1\D{x}{\uu_1}-0.225\uu_2\D{x}{\uu_2}+0.142\uu_2\D{x}{\uu_1}+0.0728\uu_1\D{x}{\uu_2}
\nonumber\\&
{}+\frac{(\uu_1-\uu_2)}{h}\left(-0.25\uu_1+0.34\uu_2\right)\D xh
\nonumber\\&
{}+\frac{1}{\re}\left(-3.84+19.3C\right)\DD{x}{\uu_1}+\frac{1}{\re}\left(2.52-6.98C\right)\DD{x}{\uu_2}\,,\label{twolayer:u1I}
\\
\mathcal L\D{t}{\uu_2}&\approx{}1.002\left(\tan\theta-\D xh\right)
\nonumber\\&
{}+\frac{1}{\re}\left(6.98\frac{\uu_1}{h^2}-5.36\frac{\uu_2}{h^2}\right)
\nonumber\\&
{}-1.25\uu_1\D{x}{\uu_1}-1.57\uu_2\D{x}{\uu_2}+0.768\uu_2\D{x}{\uu_1}+0.930\uu_1\D{x}{\uu_2}
\nonumber\\&
{}+\frac{(\uu_1-\uu_2)}{h}\D xh\left(-0.78\uu_1+0.38\uu_2\right)
\nonumber\\&
{}+\frac{1}{\re}\left(-1.98-6.98C\right)\DD{x}{\uu_1}+\frac{1}{\re}\left(5.23+5.36C\right)\DD{x}{\uu_2}\,.\label{twolayer:u2I}
\end{align}
\end{subequations}

Now we linearise this system~\eqref{eqs:twolayerI} and derive the characteristic equation, parametrised by the regularising~\(C\),
\begin{align}&
\lambda^3+\lambda^2\left[1.114\tan\theta\re ik+\frac{1}{\re}\left(1.39k^2+24.66\right)+24.66\frac{1}{\re}Ck^2\right]
\nonumber\\&\quad
-\lambda\left[0.476\tan^2\theta\re^2k^2+\tan\theta k(-2.265ik^2-14.03i-0.0144\re k)\right.
\nonumber\\&\quad
\left.+15.09\frac{1}{\re^2}k^4-0.914k^2-84.13\frac{1}{\re^2}k^2-54.73\frac{1}{\re^2}\right.
\nonumber\\&\quad
\left.-14.03C\tan\theta ik^3-\frac{1}{\re^2}(54.73C^2k^4+84.13Ck^4+109.46Ck^2)\right]
\nonumber\\&\quad
-\left[0.0725\tan^3\theta\re^3ik^3+\tan^2\theta\re k^2(-0.0191\re ik+0.615k^2+1.908)\right.
\nonumber\\&\quad
\left.+\tan\theta\frac{1}{\re}k(5.19ik^4-0.682\re^2ik^2-28.94ik^2-18.83i-0.841\re k^3\right.
\nonumber\\&\quad
\left.-0.198\re k)+\frac{1}{\re}k^2(0.209k^2-18.26)+\frac{1}{\re}(-18.83C^2\tan\theta ik^4\right.
\nonumber\\&\quad
\left.-28.94C\tan\theta ik^4-37.65C\tan\theta ik^2-18.26C k^3)+1.908C\tan^2\theta\re k^3\right.
\nonumber\\&\quad
\left.-0.198C\tan\theta k^3\right]=0\,.
\label{twolayer:eigfun3}
\end{align}

\begin{figure}
\centering
\begin{tabular}{c@{}c}
\rotatebox{90}{\hspace{10ex}growth rate~$\Re(\lambda/\re)$} &
\includegraphics{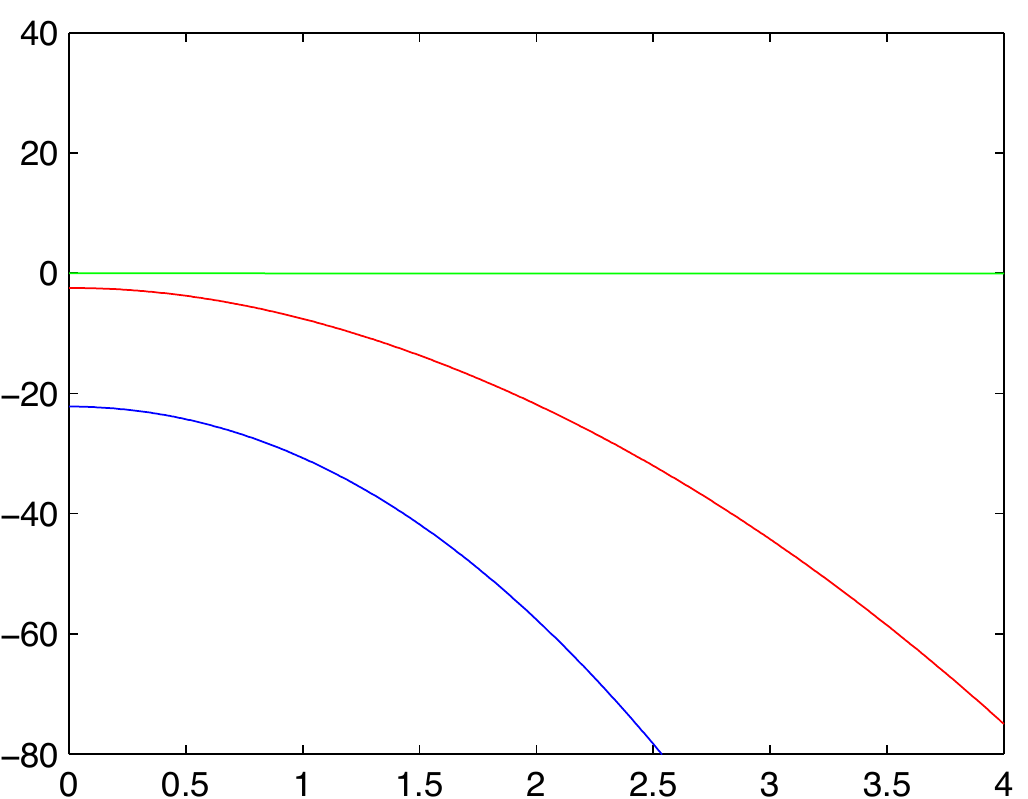}\\
&  wavenumber~$k$
\end{tabular}
\caption{Plot of the growth rate~$\Re\lambda$ varying with the nondimensional wavenumber~$k$.
The coefficient~$C=0.5$ in the regularising operator~$\mathcal L$.
The slope of the plate is~$\tan\theta=0$.
The Reynolds number is~$\re=1$.}
\label{twolayerEig3}
\end{figure}

Figure~\ref{twolayerEig3} plots the growth rate~$\Re\lambda$ varying with the nondimensional wavenumber~$k$ from the characteristic equation~\eqref{twolayer:eigfun3} for parameter $C=0.5$ in the regularising operator~$\mathcal L$.
That there is no positive growth rate demonstrates that no instability occurs, even for the high wavenumber.
The decay rates represented by the blue and red curves grow quickly for large wavenumber. 
Numerical checks indicate that the regularising coefficient $C>0.17$ to eliminate the unphysical instabilities.
This method is flexible through the wide range of choice of the positive coefficient~$C$ in the regularising operator~$\mathcal L$.
No high order derivatives, such as~$h^2\DDDD{x}{\uu_1}$, appear with this regularising method.
Thus, introducing the regularising operator~$\mathcal L$  usefully stabilises the model, and hereafter we implement the numerical simulations of the model~\eqref{eqs:twolayerI} with the regularising coefficient $C=0.5$.

\section{Gap-tooth simulation of the two layer thin fluid flow}
\label{sec:patch}

This section focuses on implementing the gap-tooth simulation of the thin fluid flow.
This section uses the two-layer model~\eqref{eqs:twolayerI} as the microscale simulator within patches.
Coupling conditions~\eqref{twolayer:u1UJ}--\eqref{twolayer:u2UJt} and~\eqref{twolayer:cpLh} are developed to couple patches together.

\begin{figure}
\begin{center}
\setlength{\unitlength}{0.8ex}
\begin{picture}(88,44)(-2,0)
\put(0,27){
\put(-2,3){\rotatebox{90}{odd \(j\)}}
\put(2,3){\vector(1,0){82}}
\put(84.5,2.5){$x$}
\setcounter{i}{-4}
\multiput(5,2.5)(12,0){7}{\line(0,1){1}%
  \stepcounter{i}%
  \put(-1,-2){$x_{j,\arabic{i}}
  \ifnum\arabic{i}=0{=}X_j\fi$}
  }
\put(23,2){\vector(-1,0){6}
  \put(6,0){\vector(1,0){6}}
  \put(7,-2){$d$}
  }
\setcounter{i}{-4}
\multiput(5,7)(12,0){7}{%
  \stepcounter{i}%
  \ifodd\arabic{i}\color{magenta}%
    \put(0,2){%
    \circle*{1}$\uu^1_{j,\arabic{i}}$
    }
     \put(0,7){%
    \circle*{1}$\uu^2_{j,\arabic{i}}$
    }
  \else\color{blue}%
    \put(0,-2){%
    \circle*{1}$h_{j,\arabic{i}}$
    }
  \fi
}
}
\put(-2,3){\rotatebox{90}{even \(j\)}}
\put(2,3){\vector(1,0){82}}
\put(84.5,2.5){$x$}
\setcounter{i}{-4}
\multiput(5,2.5)(12,0){7}{\line(0,1){1}%
  \stepcounter{i}%
  \put(-1,-2){$x_{j,\arabic{i}}
  \ifnum\arabic{i}=0{=}X_j\fi$}
  }
\put(23,2){\vector(-1,0){6}
  \put(6,0){\vector(1,0){6}}
  \put(7,-2){$d$}
  }
\setcounter{i}{-4}
\multiput(5,7)(12,0){7}{%
  \stepcounter{i}%
  \ifodd\arabic{i}\color{blue}%
    \put(0,-2){%
    \circle*{1}$h_{j,\arabic{i}}$
    }
  \else\color{magenta}%
    \put(0,2){%
    \circle*{1}$\uu^1_{j,\arabic{i}}$
    }
    \put(0,7){%
    \circle*{1}$\uu^2_{j,\arabic{i}}$
    }
  \fi
}
\end{picture}
\end{center}
\caption{Scheme of the staggered grid points of the depth $h_{j,i}$~(blue points) and the mean velocities $\uu^1_{j,i}$ and~$\uu^2_{j,i}$~(magenta points) at the \(i\)th~micro-grid point on the odd $j$th patch~(top) and the even $j$th patch~(bottom).  This diagram show the cases for \(n=5\) interior grid points in each patch (\(n'=3\)).}
\label{twolayerMicro}
\end{figure}
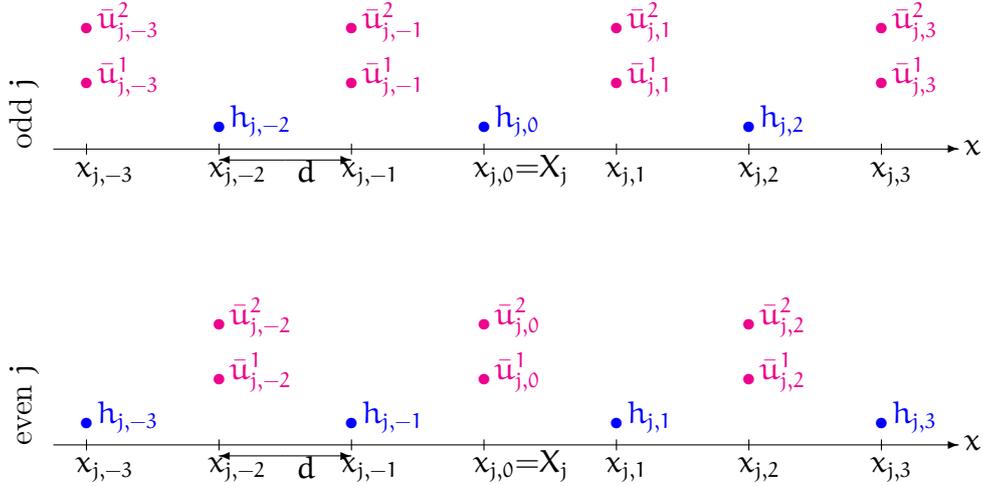

Let us focus on one patch, the~$j$th patch.
Figure~\ref{twolayerMicro} shows the staggered  grid for the depth~$h_{j,i}$~(blue points) and the layer mean velocities $\uu^1_{j,i}$ and~$\uu^2_{j,i}$~(magenta points) at the \(i\)th~micro-grid point (\(i=-n',\ldots,n'\) for \(n':=(n+1)/2\)) on the~$j$th patch.
The superscripts correspond to the two mean velocities~$\uu_1$ and~$\uu_2$ on a patch.
Let each of $m$~patches be centred on equi-spaced macroscale grid points $x=X_j=jD$, where $D=L/m$ is the macroscale spacing and $L$~is the length of the whole domain.
Each patch has relatively small width~$\lp $. 
Assume each patch has a total of~$n$ microscale interior grid points, excluding the two edge grid points, so the microscale spatial step~$d=\lp/(n+1)$.
Let each patch around a macroscale grid point~$X_j$ execute a microscale simulator.

The microscale simulator~\eqref{eqs:twolayerI} is a straightforward second order differential equations.
Approximate the \pde{}s~\eqref{eqs:twolayerI} on the $j$th patch with centred differences in microscale space of~$d$. shown in Figure~\ref{twolayerMicro}.
Use Matlab \verb|ode15s| for continuous time integration.
The dominant error comes from the macroscale coupling between patches, not the microscale discretisation~\cite[e.g.]{Cao2014}.

Define the mid patch point~$x_{j,0}=X_j$.
Let the macroscale value~$H_j=h_j(X_j,t)$ for the odd~$j$ and~$U_j=(\uu^1_j(X_j,t)+\uu^2_j(X_j,t))/2$ for even~$j$.
The inter-patch coupling uses only these macroscale values.

\subsection{Coupling conditions on the odd patches}
\label{twolayer:patchCO}

This subsection develops the coupling conditions on the odd patches.
The velocity values~$\uu^1_{j,\pm n'}$ and~$\uu^2_{j,\pm n'}$ at the edges of the odd patches are approximated by interpolation of $U_{j\pm1},U_{j\pm3},\ldots$ from neighbouring patches, which is `lifted' \cite[e.g.]{Kevrekidis09a} by requiring that the microscale dynamics lie on the slow manifold.

The dynamics in the interior of each patch is given by the microscale simulator~\eqref{eqs:twolayerI}.
The regularising operator~$\mathcal L=1-C\partial_x(h^2\partial_x)$ of the two layer mean velocities~$\uu_1(x,t)$ and~$\uu_2(x,t)$ on the odd~$j$th patch requires more details.
For the example of \(n=5\)\,, Figure~\ref{twolayerMicro}, for odd~$j$ the two momentum equations~\eqref{twolayer:u1I}--\eqref{twolayer:u2I} are of the form
\begin{align}
\mathcal L
\begin{bmatrix}\D{t}{\uu^{\ell}_{j,-3}} \\ \D{t}{\uu^{\ell}_{j,-1}} \\ \D{t}{\uu^{\ell}_{j,1}} \\ \D{t}{\uu^{\ell}_{j,3}}\end{bmatrix}
=\textsc{rhs}\,,
\label{twolayer:Lu}
\end{align}
where the~\textsc{rhs} refers to a finite difference discretisation of the right-hand sides of the two-layer \pde{}s~\eqref{eqs:twolayerI}, and the superscript~$\ell=1,2$ for the lower and upper layer.
The regularising operator~$\mathcal L$ has the matrix form, for the example \(n=5\) case, of
\begin{align}
\mathcal L=
\frac{1}{4d^2}
\begin{bmatrix} 
-Ch^2_{j,-2} &1+C(h^2_{j,-2}+h^2_{j,0}) & -Ch^2_{j,0} & 0 \\
0 & -Ch^2_{j,0}  &1+C(h^2_{j,0}+h^2_{j,2}) & -Ch^2_{j,2}
\end{bmatrix} ,
\label{twolayer:LM}
\end{align}
where~$d$ is the microscale spacial step.
The values of~$\uu^{\ell}_{j,\pm n'}$ are at the edges of the~$j$th patch.
Recall that the macroscale mean velocity $U_j(t)=(\uu^1_j(X_j,t)+\uu^2_j(X_j,t))/2$ for the even~$j$ are known.
Thus, this section completes the set of equations by finding the unknown microscale values~$\D{t}{\uu^{\ell}_{j,\pm n'}}$ and~$\uu^{\ell}_{j,\pm n'}$ from the known macroscale values~$U_j$.

The challenge is to deduce microscale values appropriate to the macroscale structures.
\cite{E:2003fk,E:2007uq} and  \cite{Malecha:2013kx} studied a heterogeneous multiscale method~(\textsc{hmm}).
The \textsc{hmm} contains two main components: an overall macroscale scheme for the macroscale variables and estimating the missing macroscale data by the microscale model.
A compression operator~$Q$ and a reconstruction operator~$R$ are defined to satisfy~$Q\uu_j=U_j$,~$RU_j=\uu_j$ and~$QRU_j=U_j$, where~$\uu_j$ is the microscale variable and~$U_j$ the macroscale variable.
Such compression and reconstruction operators combine the microscale and macroscale variables.
\cite{Kevrekidis:2003fk,Samaey:2005fk} and~\cite{Samaey2009} defined a coarse time-stepper by introducing a~\emph{lifting operator} and a corresponding~\emph{restriction operator} which transform between the microscale and macroscale variables.
These works provide the methods to relate the microscale and macroscale.

This section analogously constructs a lifting operator to give the patch edge values of~$\uu^{\ell}_{j,\pm n'}$ and~$\D{t}{\uu^{\ell}_{j,\pm n'}}$ in terms of the macroscale mean velocity~$U_j$ by assuming the system lies on the one-layer slow manifold.
Recall that the one-layer mean velocity~$\uu(X_j,t)=U_j=(\uu^1_{j,0}+\uu^2_{j,0})/2$.
One constraint on the lifting is the coupling conditions~\cite[eq.~(11)]{Cao2013} that give a linear, cubic or quintic approximation for the one layer mean velocity, such as the quintic
\begin{align}
\uu(X_j\pm rD,t)&={}\frac{1}{2}\left(U_{j+1}+U_{j-1}\right)
\pm\frac{r}{2}\left(U_{j+1}-U_{j-1}\right)
\nonumber\\&
{}+\frac{1}{16}(-1+r^2)\left(U_{j+2}-U_{j+1}-U_{j-1}+U_{j-2}\right)
\nonumber\\&
{}\pm\frac{1}{48}(-r+r^3)\left(U_{j+2}-3U_{j+1}+3U_{j-1}-U_{j-2}\right)\,.
\label{twolayer:cpL}
\end{align}
The other requirement is that the patch be on the slow manifold of macroscale waves.
The following section~\ref{twolayer:Lmodel1} details the slow manifold of the two-layer model \eqref{twolayer:u1I}--\eqref{twolayer:u2I} in term of the mean velocity~$\uu$. 
For example and for simplicity, truncate the slow manifold decription~\eqref{eqs:twolayeruU} and~\eqref{twolayer:u1Ut}--\eqref{twolayer:u2Ut} to errors~$\mathcal O(\epsilon^2)$ and obtain  the values of~$\uu^{\ell}_{j,\pm n'}$ and~$\D{t}{\uu^{\ell}_{j,\pm n'}}$ on the odd~$j$th patch as 
\begin{subequations}\label{eqs:twolayeruUJ}%
\begin{align}
\uu^1_{j,\pm n'}&= 0.587\uu(X_j\pm rD,t)+0.0129\re\tan\theta\,,\label{twolayer:u1UJ}
\\
\uu^2_{j,\pm n'}&= 1.413\uu(X_j\pm rD,t)-0.0129\re\tan\theta\,,\label{twolayer:u2UJ}
\\
\D{t}{\uu^1_{j,\pm n'}}&= -1.482\frac{1}{\re}\frac{1}{h^2}\uu(X_j\pm rD,t)+0.489\tan\theta\,,\label{twolayer:u1UJt}
\\
\D{t}{\uu^2_{j,\pm n'}}&= -3.526\frac{1}{\re}\frac{1}{h^2}\uu(X_j\pm rD,t)+1.168\tan\theta\,,\label{twolayer:u2UJt}
\end{align}
\end{subequations}
where as a leading approximation we neglect the derivatives $\D x{}\approx 0$ in the slow manifold~\eqref{eqs:twolayeruUt}.
Thus, equations~\eqref{twolayer:u1UJ}--\eqref{twolayer:u2UJ} are the coupling conditions on the odd patches.

\subsection{Coupling conditions on the even patches}
\label{twolayer:patchC}

This subsection develops the coupling conditions on the patches with even~$j$.
The values~$h_{j,\pm n'}$ at the edges of the even patches are approximated from the neighbouring macroscale grid values to give coupling conditions.

Simulate on each patch by the discretisation of the microscale \pde{}s~\eqref{eqs:twolayerI}.
The regularising operator~$\mathcal L=1-C\partial_x(h^2\partial_x)$ of the two layer mean velocities~$\uu_1(x,t)$ and~$\uu_2(x,t)$ on the even~$j$th patch requires more details.
For example, we record here details for $n=5$ as in Figure~\ref{twolayerMicro}; other~$n$ is a direct generalisation.
We need to approximate the second spatial derivatives~$\partial^2_x$ in the left-hand and right-hand sides of the \eqref{twolayer:u1I}--\eqref{twolayer:u2I} at the positions~$x_{j,\pm2}$, which requires two virtual grid values~$\uu^{\ell}_{j,\pm4}$.
We set the values $\uu^{\ell}_{j,\pm4}=\uu^{\ell}_{j,\pm2}$ on the slow manifold.
Thus, in Figure~\ref{twolayerMicro}, for even~$j$ the momentum equations~\eqref{twolayer:u1I}--\eqref{twolayer:u2I} are of the form
\begin{align}
\mathcal L
\begin{bmatrix}\D{t}{\uu^{\ell}_{j,-4}} \\ \D{t}{\uu^{\ell}_{j,-2}} \\ \D{t}{\uu^{\ell}_{j,0}}\\ \D{t}{\uu^{\ell}_{j,2}} \\ \D{t}{\uu^{\ell}_{j,4}}\end{bmatrix}
=\textsc{rhs}\,,
\label{twolayer:LuE}
\end{align}
where the regularising operator~$\mathcal L$ is discretised in the matrix form
\begin{align}
\frac{1}{4d^2}
\begin{bmatrix} 
\parbox{0.15\linewidth}{\raggedright $-Ch^2_{j,-3}$}&
\parbox{0.15\linewidth}{\raggedright $1+C(h^2_{j,-3}+h^2_{j,-1})$}& 
\parbox{0.15\linewidth}{\raggedright $-Ch^2_{j,-1}$}&
0 &
0
\\
\\
0&
\parbox{0.15\linewidth}{\raggedright $-Ch^2_{j,-1}$}&
\parbox{0.15\linewidth}{\raggedright $1+C(h^2_{j,-1}+h^2_{j,1})$}&
\parbox{0.15\linewidth}{\raggedright $-Ch^2_{j,1}$} &
0
\\
\\
0&
 0 &
\parbox{0.15\linewidth}{\raggedright $ -Ch^2_{j,1}$}&
\parbox{0.15\linewidth}{\raggedright $1+C(h^2_{j,3}+h^2_{j,1})$}&
\parbox{0.15\linewidth}{\raggedright $ -Ch^2_{j,3}$}  
\end{bmatrix} ,
\label{twolayer:LM}
\end{align}
and where the~\textsc{rhs} refers to a discretisation of the right-hand sides of the momentum equations \eqref{twolayer:u1I}--\eqref{twolayer:u2I}, the superscript~$\ell=1,2$ for the lower and upper layer, and~$d$ is the microscale spatial step.

The regularising operator~$\mathcal L$ and the \textsc{rhs}s need the values of~$h_{j,\pm n'}$. 
Coupling conditions approximate the values of~$h_{j,\pm n'}$ by interpolating the neighbouring macroscale values of~$H_{j\pm1},H_{j\pm2},\ldots$.
The coupling conditions~\cite[eq.~(11)]{Cao2013} give a linear, cubic or quintic interpolation, such as the quintic
\begin{align}
h_{j,\pm n'}&={}h(X_j\pm rD,t)=\frac{1}{2}\left(H_{j+1}+H_{j-1}\right)
\pm\frac{r}{2}\left(H_{j+1}-H_{j-1}\right)
\nonumber\\&
{}+\frac{1}{16}(-1+r^2)\left(H_{j+2}-H_{j+1}-H_{j-1}+H_{j-2}\right)
\nonumber\\&
{}\pm\frac{1}{48}(-r+r^3)\left(H_{j+2}-3H_{j+1}+3H_{j-1}-H_{j-2}\right)\,,
\label{twolayer:cpLh}
\end{align}
where~$r$ is the ratio of between the macroscale step and half of the width of a patch.
Thus the interpolation~\eqref{twolayer:cpLh}, together with equations~\eqref{twolayer:u1UJ}--\eqref{twolayer:u2UJ}, couple the patches together over the macroscale domain.

\subsection{The low order model of one layer flow}
\label{twolayer:Lmodel1}

This section derives a slow manifold of the two-layer model~\eqref{twolayer:u1I}--\eqref{twolayer:u2I} in terms of the mean velocity~$\uu$.
This slow manifold is used by sections~\ref{twolayer:patchCO}--\ref{twolayer:patchC} in order to lift the macroscale information to the microscale simulation on patches.

A fluid film model expressed in terms of the dynamics
of \emph{both} the fluid layer thickness and an overall lateral velocity (or momentum flux) is needed to resolve wave-like
dynamics in many situations \cite[]{Roberts99b}: 
falling films \cite[p.110]{Nguyen00, Chang94}; 
wave transitions \cite[]{Chang02} to solitary waves \cite[]{RuyerQuil00}; 
higher Reynolds number flows \cite[Eqn.(19)]{Prokopiou91b}; 
in rising film flow and a slot coater \cite[Eqn.(37)]{Kheshgi89}; rivulets under a sloping cylinder \cite[]{Alekseenko96}.
A slow manifold model of our two layer model corresponds to these earlier models of the fluid dynamics.

Indeed our construction here generates a slow manifold, one-layer model of the fluid dynamics which is the same as that of \cite{Roberts1996} (to the order of error of the analysis).
The distinction is that here it is derived from the two-layer model, rather than the original fluid equations: this transitivity of modelling additionally validates the modelling process.
For simplicity, we base the analysis upon a symmetric linear operator with slow eigenspace where the two layer velocities are in the ratio~$1:2$.
To derive the slow manifold of the two-layer model, we embed the physical model~\eqref{twolayer:u1I}--\eqref{twolayer:u2I} into a family of artificial problems by using another artificial parameter~$\gamma'$,
\begin{subequations}\label{eqs:twolayeruO}%
\begin{align}&
-\D{t}{\uu_1}+\textsc{rhs}_{\eqref{twolayer:u1f}}+(1-\gamma')\frac{\Eu }{\re}\frac{1}{h^2}\left(-4\uu_1+2\uu_2\right)=0\,,
\label{twolayer:u1O}
\\&
-\D{t}{\uu_2}+\textsc{rhs}_{\eqref{twolayer:u2f}}+(1-\gamma')\frac{\Eu }{\re}\frac{1}{h^2}\left(2\uu_1-\uu_2\right)=0\,,
\label{twolayer:u2O}
\end{align}
\end{subequations}
where~$\textsc{rhs}_{\eqref{twolayer:u1f}}$ and~$\textsc{rhs}_{\eqref{twolayer:u2f}}$ are the right-hand sides of the \pde{}s~\eqref{twolayer:u1f}--\eqref{twolayer:u2f}.
The variable~$\Eu $ denotes an artificial Euler parameter used to enhance convergence \cite[e.g.]{vanDyke64,Vandyke84}: computational experiments indicate that $\Eu =9/2$ delivers good convergence in the artificial parameter~$\gamma'$.
When the parameter~$\gamma'=1$, the \pde{}s~\eqref{eqs:twolayeruO} recover the original \pde{}s~\eqref{twolayer:u1f}--\eqref{twolayer:u2f}.
When the artificial parameter~$\gamma'=0$, the linear operator in the \pde{}s~\eqref{eqs:twolayeruO} is
\begin{equation*}
\frac{\Eu }{h^2\re}\begin{bmatrix} -4 &2 \\ 2& -1 \end{bmatrix},
\end{equation*}
and it is this linear operator that guides effective recursive improvements of approximations to the original \pde{}s.

For a specific model we choose to truncate to errors~$\mathcal O(\uu^4+\partial_x^4+\theta^4,{\gamma'}^7)$.
Then executing the computer algebra in Appendix~\ref{sec:app} and evaluating at the artificial parameter~$\gamma'=1$, leads to the slow manifold where the two layer velocities are
\begin{subequations}\label{eqs:twolayeruU}%
\begin{align}
\uu_1&\approx{}0.587\uu+0.0129\re h^2 \left(\tan\theta-\D xh\right)
\nonumber\\&
\quad{}-0.0468h^2\DD{x}{\uu}-0.205h\D xh\D{x}{\uu}+0.0700h\uu\DD xh
\nonumber\\&
\quad{}+\re\left(0.00465h^2\uu\D{x}{\uu}-0.0115h\uu^2\D{x}{h}-0.0105h^3\D xh\DD xh\right)\,,\label{twolayer:u1U}
\\
\uu_2&\approx{}1.413\uu-0.0129\re h^2\left(\tan\theta-\D xh\right)
\nonumber\\&
\quad{}+0.0468h^2\DD{x}{\uu}+0.205h\D xh\D{x}{\uu}-0.0700h\uu\DD xh
\nonumber\\&
\quad{}-\re\left(0.00465h^2\uu\D{x}{\uu}-0.0115h\uu^2\D{x}{h}-0.0105h^3\D xh\DD xh\right)\,.\label{twolayer:u2U}
\end{align}
\end{subequations}
Equations~\eqref{eqs:twolayeruU} are expressed in terms of the single layer depth-averaged velocity~$\uu(x,t)$ and the water depth~\(h(x,t)\).

On the slow manifold~\eqref{eqs:twolayeruU},  conservation of mass requires
\begin{equation}
\D th=-\D{x}{h\uu}\,.
\label{twolayer:hUt}
\end{equation}
The computer algebra in Appendix~\ref{sec:app} differentiates  equations~\eqref{eqs:twolayeruU} with respect to the time~$t$, and gives the rate of change of the two layer velocities in term of the mean velocity~$\uu$ and fluid depth~$h$ as
\begin{subequations}\label{eqs:twolayeruUt}%
\begin{align}
\D{t}{\uu_1}&\approx 0.489\left(\tan\theta-\D{x}{h}\right)-1.482\frac{1}{\re}\frac{\uu}{h^2}-0.904\uu\D x\uu
\nonumber\\&
\quad{}+\frac{1}{\re}\left(2.552\DD x\uu+3.077\frac{1}{h}\D xh\D x\uu-0.650\frac{\uu}{h}\DD xh\right)
\nonumber\\&
\quad{}+\re\left[0.0167h^3\DD x\uu+0.0438h^2\D xh\DD xh+0.0359h\uu\left(\D xh\right)^2\right]
\nonumber\\&
\quad{}-\re \tan\theta\left(0.0298h\uu\D{x}{h}+0.0184h^2\D{x}{\uu}\right)\,,\label{twolayer:u1Ut}
\\
\D{t}{\uu_2}&\approx 1.168\left(\tan\theta-\D{x}{h}\right)-3.526\frac{1}{\re}\frac{\uu}{h^2}-2.107\uu\D x\uu
\nonumber\\&
\quad{}+\frac{1}{\re}\left(5.701\DD x\uu+6.962\frac{1}{h}\D xh\D x\uu-0.312\frac{\uu}{h}\DD xh\right)
\nonumber\\&
\quad{}-\re\left[0.00819h^3\DD x\uu+0.0482h^2\D xh\DD xh+0.0875h\uu\left(\D xh\right)^2\right]
\nonumber\\&
\quad{}+\re \tan\theta\left(0.0847h\uu\D{x}{h}+0.0355h^2\D{x}{\uu}\right)\,.\label{twolayer:u2Ut}
\end{align}
\end{subequations}
Equations~\eqref{eqs:twolayeruUt} provide needed microscale information from the macroscale fields~$\uu$ and~$h$.

The computer algebra in Appendix~\ref{sec:app} simultaneously determines the evolution on the slow manifold of the one-layer model. 
The momentum equation is
\begin{align}
\D{t}{\uu}\approx
0.829\left(\tan\theta-\D xh\right)
-2.504\frac{1}{\re}\frac{\uu}{h^2}
-1.505\uu\D{x}{\uu}-0.151\frac{\uu^2}{h}\D xh\,.
\label{twolayer:Ut}
\end{align}
This \pde~\eqref{twolayer:Ut} has the same terms as the low order model of \cite{Roberts1996} [eq.~(11)].
The coefficients of these terms are the same to a relative error of less than~$2.2\%$.
This agreement partially verifies that the two-layer model is a reasonable model of the fluid film flow.

\section{Numerical gap-tooth simulations of the two layer thin fluid flow}
\label{sec:num}

This section explores the numerical gap-tooth simulation with the two-layer microscale simulator~\eqref{eqs:twolayerI} and the patch coupling conditions~\eqref{twolayer:cpLh} and~\eqref{twolayer:cpL}--\eqref{eqs:twolayeruUJ}.
Numerical eigenvalues and simulations show that the gap-tooth scheme reasonably simulates the macroscopic dynamics of the fluid film flow. 

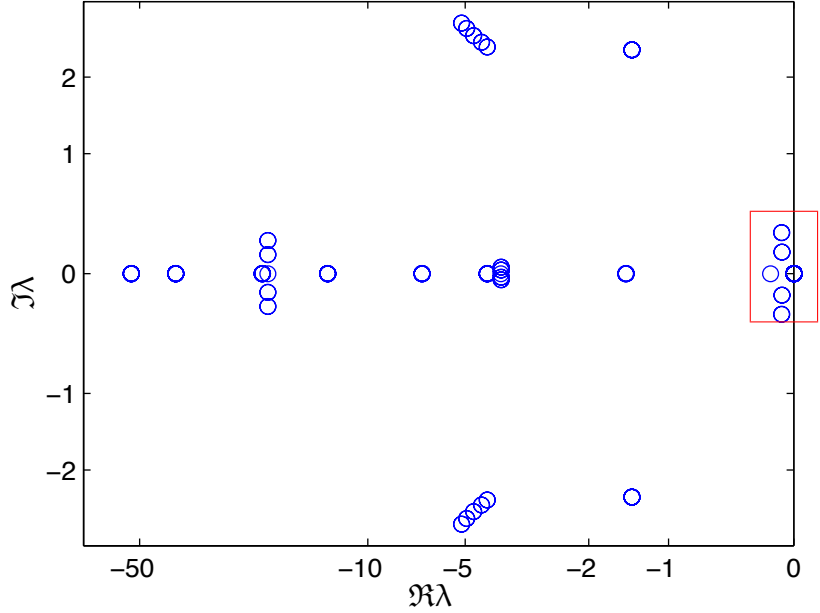
\begin{figure}
\begin{center}
\setlength{\unitlength}{0.8ex}
\begin{picture}(77,56)
\put(7,2){\includegraphics[height=54\unitlength]{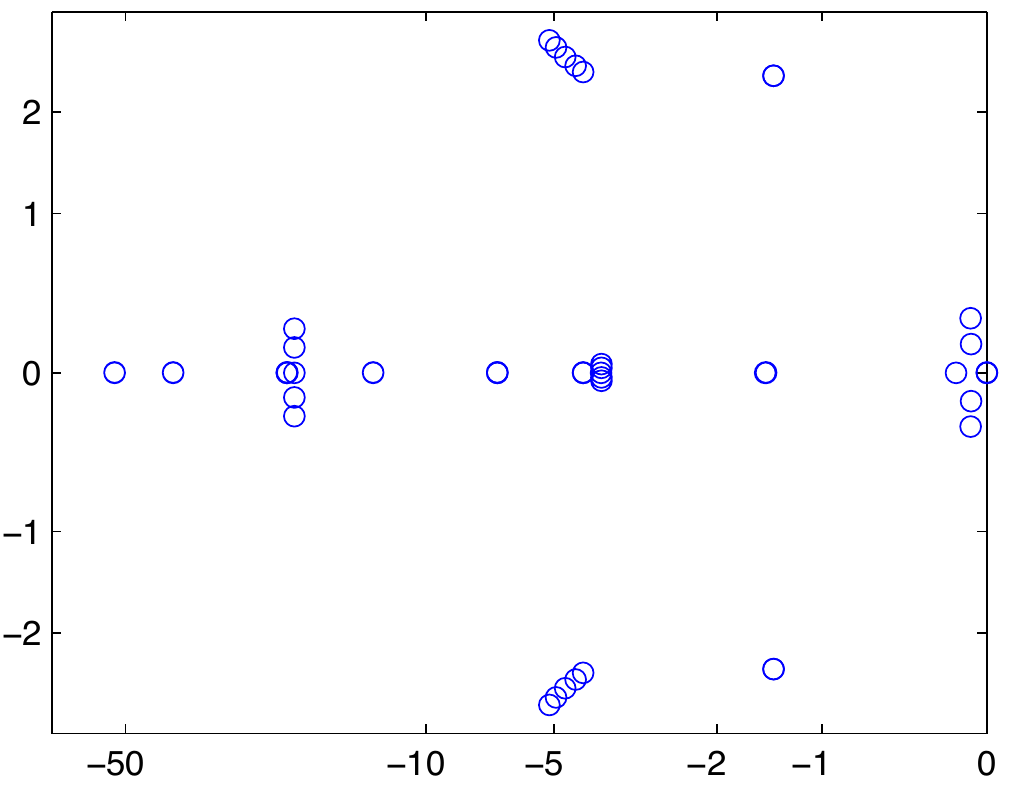}}
\put(4,2){\rotatebox{90}{\hspace{20ex}$\Im\lambda$}}
\put(71,26){\color{red}\framebox(6,10)}
\put(40,0){$\Re\lambda$}
\end{picture}
\end{center}
\caption{Plots of the growth rate~$\Re\lambda$ versus the frequency~$\Im\lambda$ in the gap-tooth simulation of the thin film flow (with nonlinear axis scaling) over horizontal plate.
There are $m=10$ patches and $n=9$ microscale grid points within a patch.
The Reynolds number $\re=15$ and the coefficient in the regularising operator~$\mathcal L$ is $C=0.5$.
}
\label{twolayerPeig}
\end{figure}%

\begin{figure}
\centering
\begin{tabular}{c@{}c}
\rotatebox{90}{\hspace{20ex}$\Im\lambda$} &
\includegraphics{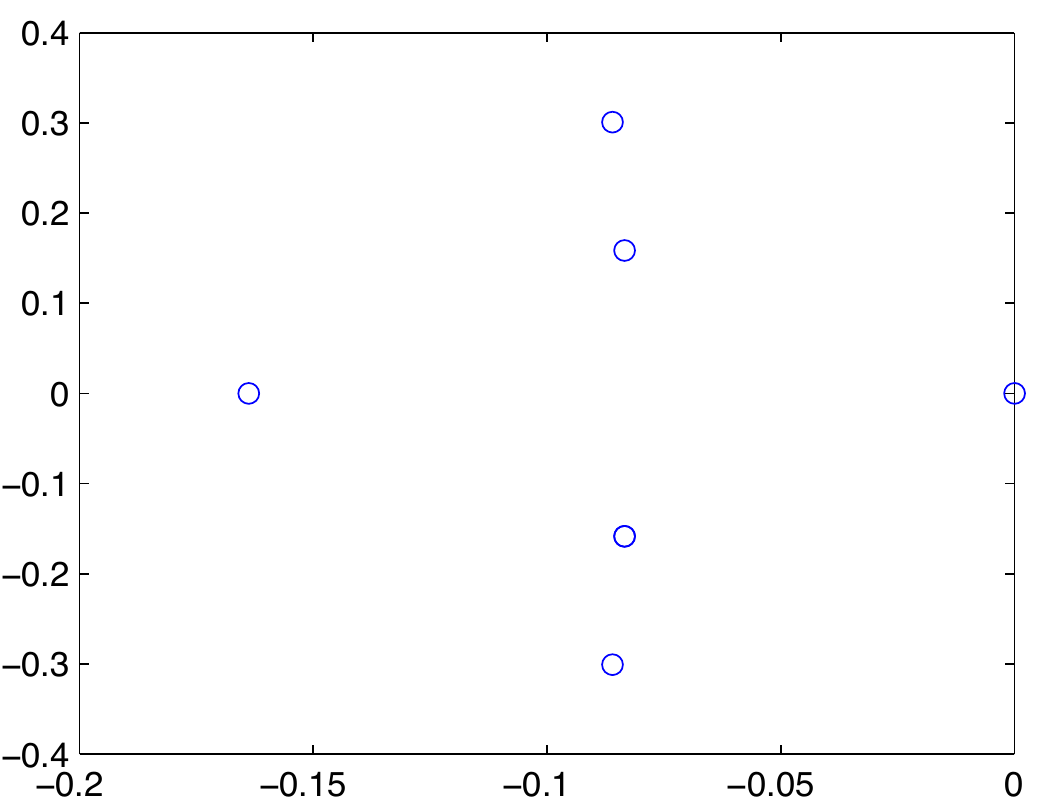}\\
& $\Re\lambda$
\end{tabular}
\caption{Zoom into the values in the red rectangle in Figure~\ref{twolayerPeig}.  These eigenvalues which represent the macroscale modes.}
\label{twolayerPeigMa}
\end{figure}%

Consider the thin fluid flow on a horizontal plate: that is, the mean slope $\tan\theta=0$\,.
Distribute~$m$ patches in the macroscale domain of length $L=m\pi$, so the distance between the neighbouring patches is $D=L/m=\pi$.
Divide each patch into~$n+1$ equal microscale intervals by~$n+2$ grid points, then the distance between neighbouring microscale points is $d=2rD/(n+1)$, where $r$~is the ratio between the half-width of each patch and the macroscale inter-patch distance~$D$.
We approximate the spatial derivatives in the right-hand side of the two-layer \pde{}s~\eqref{eqs:twolayerI} by centred differences on the staggered microscale grids.

Figure~\ref{twolayerPeig} plots the growth rate~$\Re\lambda$ versus the frequency~$\Im\lambda$ in the gap-tooth simulation of the thin fluid flow. 
There are~$m=10$ patches and~$n=9$ microscale grids on a patch.
The plate has a length~$L=m\pi\approx31.4$, so the distance between the neighbouring patches is $D=L/m=3.14$.
The ratio $r=1/6$, so the microscale step on a patch is $d=2rD/(n+1)\approx0.17$.
The Reynolds number~$\re=15$ and the coefficient in the regularising operator~$\mathcal L$ is~$C=0.5$.
The negative growth rates imply the waves decay in time.

There are $90$ pairs of eigenvalues for the system:
\begin{itemize}
\item the eigenvalues with large decay rates ($\Re\lambda<-1$) and zero imaginary parts  predominantly represent the viscous decay of the lateral shear modes;
\item the eigenvalues with large decay rates ($\Re\lambda<-1$) but with large imaginary parts ($|\Im\lambda|>2$) correspond to microscale waves within each patch; and 
\item the eigenvalues with small growth rates ($|\Re\lambda|<0.2$) in the red rectangle, zoomed in by Figure~\ref{twolayerPeigMa}, represent the interesting macroscale wave-like dynamics.  
This set of ten eigenvalues includes two groups of four eigenvalues representing waves of wavenumbers~\(0.2\) and~\(0.4\), a zero eigenvalue representing conservation of water, and $\lambda\approx-0.16$ representing the decay of homogeneous shear.
\end{itemize}
The pattern of eigenvalues seen in Figure~\ref{twolayerPeig} is typical over a wide range of parameters.
It shows the emergence of the macroscale wave modes from among the fast microscale modes within each patch.

\begin{figure}
\centering
\begin{tabular}{rc}
$t=0$\\\rotatebox{90}{\hspace{8ex}$h$} &
\includegraphics{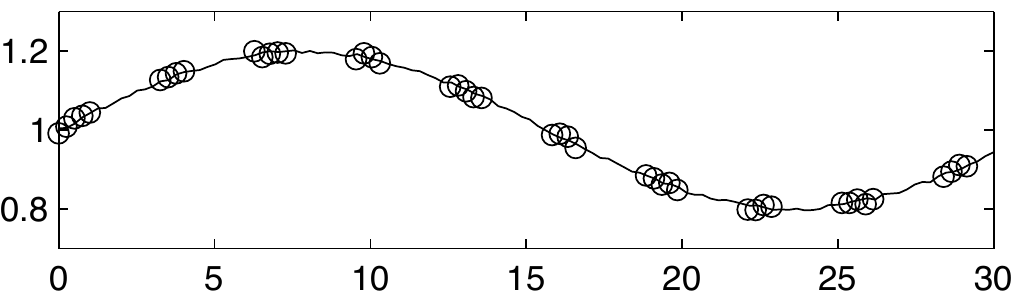}\\
& $x$
\\
$t=2$\\\rotatebox{90}{\hspace{8ex}$h$ } &
\includegraphics{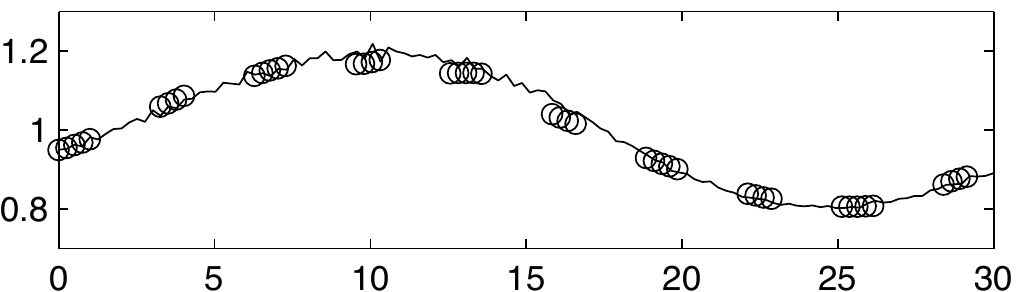}\\
& $x$
\\
$t=10$\\\rotatebox{90}{\hspace{8ex}$h$  } &
\includegraphics{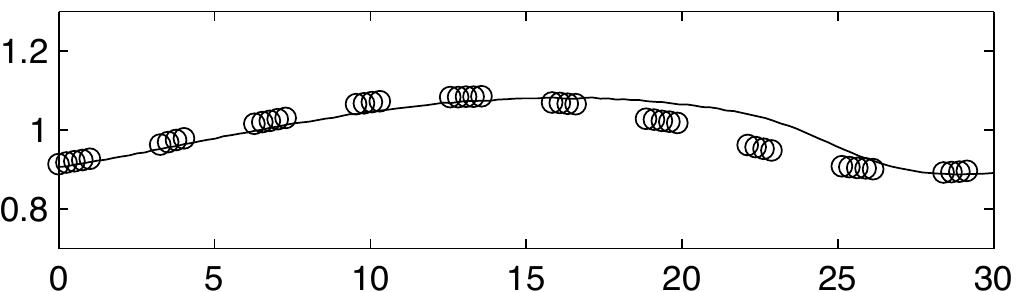}\\
& $x$
\\
$t=20$\\\rotatebox{90}{\hspace{8ex}$h$  } &
\includegraphics{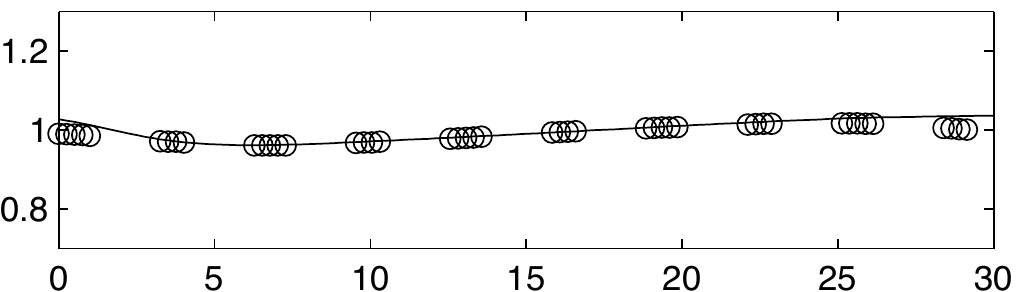}\\
& $x$
\end{tabular}
\caption{Plots of the fluid thickness of the thin fluid flow on the domain~$[0\ 10\pi]$: (circle) for the gap-tooth simulation; and (line) for microscale simulation on the whole domain.
There are $m=10$ patches and $n=9$ microscale grids on a patch.
The Reynolds number $\re=15$ and the regularising coefficient~$C=0.5$.
The plate has a slope~$\tan\theta=0$.
Figure~\ref{twolayerCu} shows the corresponding fluid velocities.
}
\label{twolayerCh}
\end{figure}

\begin{figure}
\centering
\begin{tabular}{rc}
$t=0$\\
\rotatebox{90}{\hspace{4ex}$\uu_1\ \&\ \uu_2$} &
\includegraphics{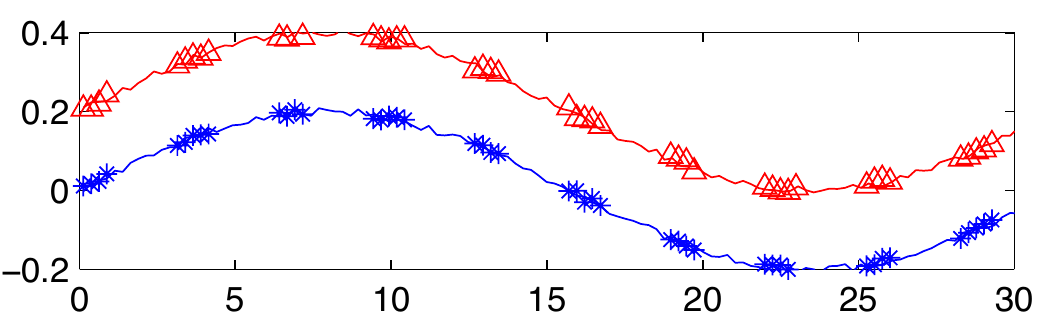}\\
& $x$
\\$t=2$\\
\rotatebox{90}{\hspace{4ex}$\uu_1\ \&\ \uu_2$ } &
\includegraphics{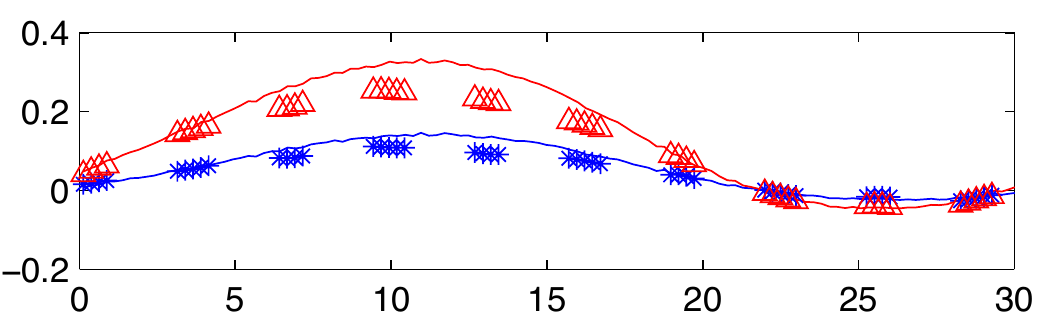}\\
& $x$
\\$t=10$\\
\rotatebox{90}{\hspace{4ex}$\uu_1\ \&\ \uu_2$  } &
\includegraphics{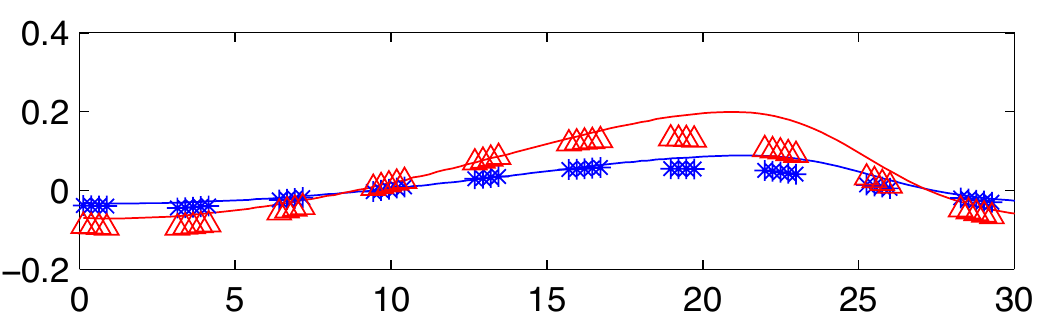}\\
& $x$
\\$t=20$\\
\rotatebox{90}{\hspace{4ex}$\uu_1\ \&\ \uu_2$ } &
\includegraphics{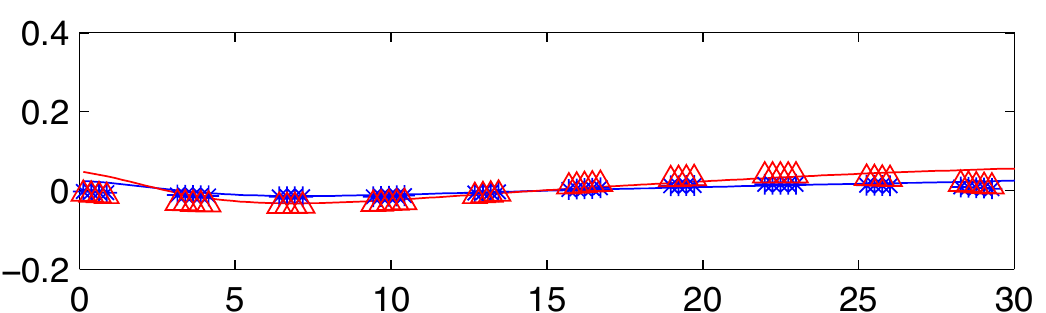}\\
& $x$
\end{tabular}
\caption{Plots of the layer mean velocities: (blue line) represents~$\uu_1$ and (red line) represents~$\uu_2$ in the microscale simulation over the whole domain; and (blue star) represents~$\uu_1$ and~(red triangle) represents~$\uu_2$ in the gap-tooth simulation of the fluid film flow on the domain~$[0\ 10\pi]$.
There are~$m=10$ patches and~$n=9$ microscale grids on a patch.
}
\label{twolayerCu}
\end{figure}

Simulations confirm the nonlinear dynamics of the gap-tooth scheme are appropriate.
Figure~\ref{twolayerCh} plots and compares the free surface of the thin film flow in a gap-tooth simulation and in the microscale simulation over the whole domain.
Time integration invoked the Matlab~\verb|ode15s| function.
There are $m=10$ patches and $n=9$ microscale grids on a patch.
The $t=0$ graph shows that initially impose a perturbation~$0.2\sin(2\pi/Lx)$ with small random noises to the equilibrium of fluid thickness one and initial layer mean velocities $\uu_1=0$ and $\uu_2=0.2$, so the mean layer velocity $\uu=(\uu_1+\uu_2)/2=0.1$.
The $t=2$ to $t=10$ graphs show that the microscale modes on a patch smooth quickly and the macroscale waves propagate and slowly decay, which agrees with the large decay rates of the microscale waves and small decay rates of the macroscale modes in Figure~\ref{twolayerPeig}.
Further numerical simulation shows that the macroscale waves decay to the equilibrium of thickness one after time $t\approx30$ for the Reynolds number $\re=15$.
When the Reynolds number~$\re$ becomes smaller, the macroscale waves decay faster to the equilibrium due to the stronger viscous effects.
Figure~\ref{twolayerCh} shows that the gap-tooth simulation agrees with the microscale simulation over the whole domain.
The likely reason for  the  differences in $t=10$ graph is that the gap-tooth simulation near the water bore involves the error $\mathcal O(D^2)$ for the macroscale step $D=L/m\approx3.14$, while the microscale simulation on the whole macroscale domain involves the error $\mathcal O(d^2)$ for the microscale step $d=2rD/(n+1)=0.13$\,: the ratio of these errors is \(\Ord{D^2/d^2}\approx 500\) showing the macroscale simulation is limited largely by the interpolation between patches.
As the number of patches increases, the error~$\mathcal O(D^2)$ declines, and the gap-tooth simulation performs better.
We show just \(m=10\) patches here for clarity.

Figure~\ref{twolayerCu} plots the corresponded layer mean velocities~$\uu_1$~(blue) and~$\uu_2$~(red) in the gap-tooth simulation and in the microscale simulation on the whole domain.

\section{Conclusion}
\label{sec:con}

This work applies the gap-tooth scheme to the thin fluid flow with a derived novel two-layer model~\eqref{eqs:twolayerI} being the microscale simulator.
The two-layer model~\eqref{eqs:twolayerI} provides more microscale modes, but without the full complexity of fully resolved vertical structures.

Based on centre manifold theory, section~\ref{sec:micro}--\ref{sec:eig} derive the two-layer model~\eqref{eqs:twolayerI} from the continuity and modified Navier--Stokes equations.
The model~\eqref{eqs:twolayerI} includes the effects of gravity, drag, advection and dispersion.
The slow manifold of the two-layer model agrees with the one-layer model of  thin film flow by~\cite{Roberts1996}.
Eigenvalues analysis in section~\ref{sec:eig} indicates unphysical instability appear, but a consistent regularising operator~$\mathcal L$ stabilises the model.

Then we simulate the thin fluid flow by new adaptations of the gap-tooth scheme with the two-layer model~\eqref{eqs:twolayerI} being the microscale simulator.
Section~\ref{twolayer:Lmodel1} algebraically derives the slow manifold of the two-layer model~\eqref{eqs:twolayerI} in order to lift the macroscale information to the microscale simulation on patches.
Then coupling conditions~\eqref{twolayer:cpLh} and~\eqref{twolayer:u1UJ}--\eqref{twolayer:u2UJ} are developed to couple the patches.
Numerical eigenvalues in Figure~\ref{twolayerPeig}--\ref{twolayerPeigMa} show that the macroscale modes decay slowly over the whole domain, while the microscale modes oscillate fast and decay quickly to quasi-equilibrium quickly.
Non-zero frequencies~$\Im\lambda$ in Figure~\ref{twolayerPeigMa} indicate waves are supported on the free surface over the whole domain.
Numerical simulations in Figure~\ref{twolayerCh}--\ref{twolayerCu} show that the gap-tooth scheme with the two-layer model~\eqref{eqs:twolayerI} being microscale simulator work well enough for the thin film flow.
Future work is planned to build on this base in order to apply the gap-tooth scheme to the complicated microscale physics of turbulent, possibly ice covered, water waves, with the aim of empowering effective computation over very large lateral scales.

\bibliographystyle{agsm}
\IfFileExists{ajr.sty}
{\bibliography{bibexport,ajr,bib}}
{\bibliography{bibexport,Turbulence}}

\appendix
\section{Ancillary computer algebra program}
\label{sec:app}

This appendix lists the computer algebra to derive a two-layer model of the viscous layer of fluid.
Denote fluid thickness~$h(x,t)$ by~\verb|h|, layer mean velocities~$\uu_j(x,t)$ by~\(\verb|uu|j\) for the lower layer \(j=1\) and upper layer \(j=2\), and their evolution $\D th=\verb|gh|$ and $\D t{\uu_j}=\verb|gu|j$.  
The Reynolds number is~\verb|re|, and the coefficients of lateral and normal gravitational forcing are \verb|grx|~and~\verb|grz:=1|.  

Use the operator \verb|h(m)| to denote $m$~lateral derivatives of the fluid thickness~$h$, $\partial_x^mh$, 
and similarly \verb|uuj(m)| denotes $m$~lateral derivatives of the layer mean velocity~$\uu_j$, $\partial_x^m\uu_j$.  
Use~\verb|d| to count the number of lateral $x$~derivatives so we can easily truncate the asymptotic expansion.
These operators depend upon time and lateral space.  
Then the spatial derivative $\partial_x\verb|h(m)|=\verb|h(m+1)|$, and the time derivative $\partial_t\verb|h(m)|= \partial_x^m\verb|gh|$, for example. 

\small
\verbatimlisting{CATwoLayer.txt}

\end{document}